\documentclass[10pt]{article}

\newlength{\minitwocolumn}
\font\teneufm=eufm10
\font\seveneufm=eufm7
\font\fiveeufm=eufm5
\newfam\eufmfam
\textfont\eufmfam=\teneufm
\scriptfont\eufmfam=\seveneufm
\scriptscriptfont\eufmfam=\fiveeufm

\makeatletter
\@addtoreset{equation}{section}
\makeatother

\newtheorem{thm}{Theorem}[section]
\newtheorem{prop}[thm]{Proposition}

\newtheorem{cor}[thm]{Corollary}

\newtheorem{dfn}[thm]{Definition}

\title{
\Large{\bf
SCREENINGS AND VERTEX OPERATORS OF\\
QUANTUM SUPERALGEBRA
$U_q(\widehat{{sl}}(N|1))$
}}
\begin{document}

\maketitle

\begin{center}
{TAKEO KOJIMA}
\\~\\
{\it
Faculty of Engineering,
Yamagata University, Jonan 4-3-16, Yonezawa 992-8510,
Japan\\
kojima@yz.yamagata-u.ac.jp}
\end{center}

~\\
~\\

\begin{abstract}
We construct the screening currents
of the quantum superalgebra $U_q(\widehat{sl}(N|1))$
for an arbitrary level $k \neq -N+1$.
We show that these screening currents commute with 
the superalgebra modulo total difference.
We propose bosonizations of the vertex operators 
by using
the screening currents.
We check that these vertex operators are the intertwiners among
the Fock-Wakimoto representation and the typical representation 
for rank $N \leq 4$.
\end{abstract}

\newpage

\section{Introduction}

Bosonizations are known to be a powerful method
to construct correlation functions 
not only in conformal field theory
\cite{Bouwknegt-McCarthy-Pilch1},
but also in exactly solvable lattice model \cite{Jimbo-Miwa}.
In the previous paper \cite{Kojima1}
we constructed a bosonization
of the quantum affine superalgebra 
$U_q(\widehat{sl}(N|1))$ for an arbitrary level 
$k \in {\bf C}$.
Bosonizations for an arbitrary level $k \in {\bf C}$
\cite{Wakimoto, Feigin-Frenkel, Matsuo,
Shiraishi, Awata-Odake-Shiraishi1, 
Awata-Odake-Shiraishi2, Zhang-Gould, Kojima1}
are completely different from those of level $k=1$
\cite{Frenkel-Kac, Segal, Frenkel-Jing,
Bernard, Jing-Koyama-Misra,
Jing1, Jing2, Kimura-Shiraishi-Uchiyama, Zhang, Yang-Zhang}.
This paper is a continuation of the paper 
\cite{Kojima1}.
In this paper we focus our attention on the screening
currents, that play an important role
in level $k$ bosonizations \cite{Kojima1} .
We construct the screening currents that commute with 
the quantum superalgebra $U_q(\widehat{sl}(N|1))$ 
modulo total difference, for an arbitrary level $k \neq -N+1$.
Using the screening currents, we construct
the screening operators 
that commute with the quantum superalgebra.
The screening currents are useful to
study level $k$ bosonizations, that isn't irreducible representation.
For instance,
(1) the screening currents balance the "background charge" of 
the vertex operators \cite{Matsuo, Zhang-Gould,
Dotsenko-Fateev, Konno, Kato-Quano-Shiraishi}, and 
(2) the irreducible representation is constructed
from the Felder complex by
the screening currents \cite{Felder, Bernard-Felder, 
Bouwknegt-McCarthy-Pilch2, Konno, Feigin-Jimbo-Miwa-Odesskii-Pugai}.  
In this paper we focus our attention on
the background charge problem.
We propose bosonizations of
the vertex operators \cite{Frenkel-Reshetikhin} 
that are the intertwiners among
the Fock-Wakimoto module and
the typical representation, by using the screening operators.
We check the intertwining property of
these bosonizations of the vertex operators for rank $N \leq 4$.
The screening currents and the vertex operators have been 
constructed only for
$U_q(\widehat{sl}(N))$, $U_q(\widehat{sl}(2|1))$ 
\cite{Matsuo, Shiraishi, Awata-Odake-Shiraishi1, Zhang-Gould} by now.
This paper gives a higher-rank generalization of
the screenings and the vertex operators
in $U_q(\widehat{sl}(2|1))$ paper \cite{Zhang-Gould}.
The representation theories of the superalgebra are much more complicated 
than non-superalgebra and have rich structures
\cite{Kac1, Kac2, Kac-Wakimoto, Kac3}.

This paper is organized as follows.
In section 2 we recall the Chevalley realization and 
the Drinfeld realization of
the quantum affine superalgebra $U_q(\widehat{sl}(N|1))$.
In section 3 we give the bosonization of
the quantum affine superalgebra $U_q(\widehat{sl}(N|1))$
for an arbitrary level $k$.
We propose the Fock-Wakimoto module by the $\xi$-$\eta$ system.
In section 4 we introduce the screening currents that
commute with the superalgebra modulo total difference,
for an arbitrary level $k \neq -N+1$.
In section 5 we propose bosonizations of
the vertex operators of
$U_q(\widehat{sl}(N|1))$.
We give the level-zero representation of the Drinfeld
generators for $U_q(\widehat{sl}(3|1))$ in this section
(resp. $U_q(\widehat{sl}(4|1))$ in appendix \ref{appendixB}).
We check that the vertex operators 
are the intertwiners among
the Fock-Wakimoto realization and the typical representation
of the quantum
superalgebra $U_q(\widehat{sl}(N|1))$ for small rank
$N \leq 4$.
We show non-vanishing property of the correlation functions. 
In appendix A we summarize useful formulae of 
the normal orderings.
In appendix B we summarize the level-zero representation
of the Drinfeld generators for $U_q(\widehat{sl}(4|1))$.

\section{Quantum affine superalgebra 
$U_q(\widehat{sl}(N|1))$}

In this section
we recall the definition of 
the quantum superalgebra $U_q(\widehat{sl}(N|1))$.
Throughout this paper we fix a complex number
$0<|q|<1$.

\subsection{Chevalley generator}

We recall the definition of the quantum superalgebra
$U_q(\widehat{sl}(N|1))$ $(N=2,3,\cdots)$ 
in terms of the Chevalley generators
\cite{Yamane}.
The Cartan matrix of the affine superalgebra
$\widehat{sl}(N|1)$ is given by

\begin{eqnarray}
(A_{i,j})_{0\leq i,j \leq N}=
\left(\begin{array}{ccccccc}
0&-1&0&\cdots&\cdots&0&1\\
-1&2&-1&\cdots&\cdots\cdots&&0\\
0&-1&2&\cdots&\cdots&\cdots&\cdots\\
\cdots&\cdots&\cdots&\cdots&\cdots&\cdots&\cdots\\
\cdots&\cdots&\cdots&\cdots&2&-1&0\\
0&\cdots&\cdots&\cdots&-1&2&-1\\
1&0&\cdots&\cdots&0&-1&0
\end{array}\right).
\end{eqnarray}

~\\
We introduce the orthonormal basis $\{\epsilon_i|
i=1,2,\cdots,N+1\}$
with the bilinear form,
$(\epsilon_i|\epsilon_j)=\nu_i \delta_{i,j}$,
where $\nu_j=+~(j=1,2,\cdots,N)$ and $\nu_{N+1}=-$.
Define 
$\bar{\epsilon}_i=\epsilon_i-\frac{\nu_i}{N-1}
\sum_{j=1}^{N+1}\epsilon_j$.
Note that $\sum_{j=1}^N \bar{\epsilon}_j=0$.
The classical simple roots $\bar{\alpha}_i$
and the classical fundamental weights 
$\bar{\Lambda}_i$
are defined by
$\bar{\alpha}_i=\nu_i \epsilon_i-\nu_{i+1} \epsilon_{i+1}$,
$\bar{\Lambda}_i=\sum_{j=1}^i \bar{\epsilon}_j$
$(1\leq i \leq N)$.
Introduce
the affine weight $\Lambda_0$
and the null root
$\delta$ satisfying
$(\Lambda_0|\Lambda_0)=(\delta|\delta)=0$,
$(\Lambda_0|\delta)=1$, $(\Lambda_0|\epsilon_i)=0$,
$(\delta|\epsilon_i)=0$, $(1\leq i \leq N)$.
The other affine weights
and the affine roots
are given by
$\alpha_0=\delta-\sum_{j=1}^N \bar{\alpha}_j$,
$\alpha_i=\bar{\alpha}_i$,
$\Lambda_i=\bar{\Lambda}_i+\Lambda_0$,
$(1\leq i \leq N)$.
Let
$P=\oplus_{j=1}^N{\bf Z}\Lambda_j
\oplus {\bf Z}\delta$
and $P^*=\oplus_{j=1}^N{\bf Z}h_j
\oplus {\bf Z}d$
the affine $\widehat{sl}(N|1)$ weight lattice 
and its dual lattice, respectively.

\begin{dfn}~\cite{Yamane}~
The quantum affine
superalgebra $U_q(\widehat{sl}(N|1))$ are
generated by
the Chevalley generators
$h_i, e_i, f_i~(1\leq i \leq N)$.
The ${\bf Z}_2$-grading
of the generators are $|e_0|=|f_0|=|e_N|=|f_N|=1$
and zero otherwise.
The defining relations are
\begin{eqnarray}
&&[h_i,h_j]=0,~~[h_i,e_j]=A_{i,j}e_j,
~~[h_i,f_j]=-A_{i,j}f_j,~~
[e_i,f_j]=\delta_{i,j}\frac{q^{h_i}-q^{-h_i}}{q-q^{-1}},
\end{eqnarray}
and the Serre relations
\begin{eqnarray}
&&~[e_j,[e_j,e_i]_{q^{-1}}]_q=0,~~
~[f_j,[f_j,f_i]_{q^{-1}}]_q=0
~~for~ |A_{i,j}|=1, i \neq 0,N.
\end{eqnarray}
Here and throughout this paper,
we use the notations
\begin{eqnarray}
~[X,Y]_\xi=XY-(-1)^{|X||Y|}\xi YX.
\end{eqnarray}
We write $[X,Y]_1$ as $[X,Y]$ for simplicity.
The quantum affine superalgebra $U_q(\widehat{sl}(N|1))$
has the ${\bf Z}_2$-graded Hopf-algebra structure.
We take the following coproduct
\begin{eqnarray}
\Delta(e_i)=e_i\otimes 1+q^{h_i}\otimes e_i,~~
\Delta(f_i)=f_i\otimes q^{-h_i}+1 \otimes f_i,~~
\Delta(h_i)=h_i \otimes 1+1 \otimes h_i, 
\end{eqnarray}
and the antipode
\begin{eqnarray}
S(e_i)=-q^{-h_i}e_i,~~
S(f_i)=-f_i q^{h_i},~~
S(h_i)=-h_i.
\end{eqnarray}
The coproduct $\Delta$ satisfies an algebra automorphism
$\Delta(XY)=\Delta(X)\Delta(Y)$
and the antipode $S$ satisfies
a ${\bf Z}_2$-graded algebra anti-automorphism
$S(XY)=(-1)^{|X||Y|}S(Y)S(X)$.
The multiplication rule 
for the tensor product is ${\bf Z}_2$-graded 
and is defined for homogeneous elements
$X,Y,X',Y' \in U_q(\widehat{sl}(N|1))$ and
$v \in V, w \in W$ by
$X \otimes Y \cdot X' \otimes Y'=(-1)^{|Y||X'|}
X X' \otimes Y Y'$ and
$X \otimes Y \cdot v \otimes w=(-1)^{|Y||v|}
X v \otimes Y w$,
which extends to inhomogeneous elements through linearity.
\end{dfn}
We sometimes use 
the anti-commutator $\{X,Y\}=XY+YX
=[X,Y]_1$ for $|X|=|Y|=1$.

\subsection{Drinfeld realization}

We recall the Drinfeld's second realization 
of the quantum affine superalgebra $U_q(\widehat{sl}(N|1))$
\cite{Yamane, Drinfeld}.
The Drinfeld realization is convenient for
constructions of bosonizations.
We use the standard symbol of $q$-integer
\begin{eqnarray}
~[a]=\frac{q^a-q^{-a}}{q-q^{-1}}.
\end{eqnarray}

\begin{dfn}~~\cite{Yamane}~
The Drinfeld generators of
the quantum affine
superalgebra $U_q(\widehat{sl}(N|1))$
are
$X_{i,m}^\pm,~h_{i,m}$, $c$  $(1\leq i \leq N,
m \in {\bf Z})$.
The ${\bf Z}_2$-grading
of the Drinfeld generators
are : $|X_{N,m}^\pm|=1$ $(m \in {\bf Z})$
and zero otherwise.
Defining relations are
\begin{eqnarray}
&&~c : {\rm central},~[h_i,h_{j,m}]=0,\\
&&~[h_{i,m},h_{j,n}]=\frac{[A_{i,j}m] [cm]}{m}
\delta_{m+n,0}~~(m,n\neq 0),\\
&&~[h_i,X_j^\pm(z)]=\pm A_{i,j}X_j^\pm(z),\\
&&~[h_{i,m}, X_j^+(z)]=\frac{[A_{i,j}m]}{m}
q^{-\frac{c}{2}|m|} z^m X_j^+(z)~~(m \neq 0),\\
&&~[h_{i,m}, X_j^-(z)]=-\frac{[A_{i,j}m]}{m}
q^{\frac{c}{2}|m|}z^m X_j^-(z)~~(m \neq 0),\\
&&(z_1-q^{\pm A_{i,j}}z_2)
X_i^\pm(z_1)X_j^\pm(z_2)
=
(q^{\pm A_{j,i}}z_1-z_2)
X_j^\pm(z_2)X_i^\pm(z_1)~~~{\rm for}~|A_{i,j}|\neq 0,
\\
&&
~[X_i^\pm(z_1),X_j^\pm(z_2)]=0~~~{\rm for}~|A_{i,j}|=0,
\label{def:Drinfeld8}\\
&&~[X_i^+(z_1),X_j^-(z_2)]
=\frac{\delta_{i,j}}{(q-q^{-1})z_1z_2}
\left(
\delta(q^{-c}z_1/z_2)\Psi_i^+(q^{\frac{c}{2}}z_2)-
\delta(q^{c}z_1/z_2)\Psi_i^-(q^{-\frac{c}{2}}z_2)
\right), \\
&& 
~\left[X_i^\pm(z_{1}),
\left[X_i^\pm(z_{2}),X_j^\pm(z)\right]_{q^{-1}}
\right]_q
+\left(z_1 \leftrightarrow z_2\right)=0
~~~{\rm for}~|A_{i,j}|=1,~i\neq N.
\end{eqnarray}
where we have used
$\delta(z)=\sum_{m \in {\bf Z}}z^m$.
Here we have 
used the abbreviation $h_i={h_{i,0}}$.
We have
used the generating function
\begin{eqnarray}
X_j^\pm(z)&=&
\sum_{m \in {\bf Z}}X_{j,m}^\pm z^{-m-1},\\
\Psi_i^+(z)&=&q^{h_i}
\exp\left(
(q-q^{-1})\sum_{m>0}h_{i,m}z^{-m}
\right),\\
\Psi_i^-(z)&=&q^{-h_i}
\exp\left(-(q-q^{-1})\sum_{m>0}h_{i,-m}z^m\right).
\end{eqnarray}
\end{dfn}

The relations between
the Chevalley generators and
the Drinfeld realization are given by
\begin{eqnarray}
&&h_{i}=h_{i,0},~~~e_i=X_{i,0}^+,~~~f_i=X_{i,0}^-
~~~{\rm for}~~1\leq i \leq N,\\
&&h_0=c-(h_{1,0}+\cdots+h_{N,0}),\\
&&e_0=(-1)[X_{N,0}^- \cdots,[
X_{3,0}^-,[X_{2,0}^-,X_{1,1}^-]_{q^{-1}}]_{q^{-1}}\cdots ]_{q^{-1}}
q^{-h_{1,0}-h_{2,0}-\cdots-h_{N,0}},\\
&&f_0=q^{h_{1,0}+h_{2,0}+\cdots+h_{N,0}}
[\cdots
[[X_{1,-1}^+,X_{2,0}^+]_q, X_{3,0}^+ ]_q, \cdots X_{N,0}^+]_q.
\end{eqnarray}

\section{Bosonization of $U_q(\widehat{sl}(N|1))$}

In this section we recall the bosonization of $U_q(\widehat{sl}(N|1))$
for an arbitrary level $k \in {\bf C}$ \cite{Kojima1}.

\subsection{Boson}

We introduce the bosons
and the zero-mode operators
$a_m^j, Q_a^j$ $(m \in {\bf Z},
1\leq j \leq N)$, 
$b_m^{i,j}, Q_b^{i,j}$,
$c_m^{i,j}, Q_c^{i,j}$
$(m \in {\bf Z}, 1\leq i<j \leq N+1)$.
The bosons $a_m^i, b_m^{i,j}, c_m^{i,j}$, 
$(m \in {\bf Z}_{\neq 0})$ 
and the zero-mode operators
$a_0^i,Q_a^i$, $b_0^{i,j},Q_b^{i,j}$,
$c_0^{i,j}, Q_c^{i,j}$ that satisfy
\begin{eqnarray}
&&~[a_m^i,a_n^j]=\frac{[(k+N-1)m] [A_{i,j}m]}{m}
\delta_{m+n,0},~~
[a_0^i, Q_a^j]=(k+N-1)A_{i,j},
\\
&&~[b_m^{i,j},b_n^{i',j'}]=
-\nu_i \nu_j \frac{[m]^2}{m}
\delta_{i,i'}\delta_{j,j'}\delta_{m+n,0},~~
[b_0^{i,j},Q_b^{i',j'}]=
-\nu_i \nu_j \delta_{i,i'}\delta_{j,j'},
\\
&&~[c_m^{i,j},c_n^{i',j'}]=
\nu_i \nu_j \frac{[m]^2}{m}
\delta_{i,i'}\delta_{j,j'}
\delta_{m+n,0},~~
[c_0^{i,j},Q_c^{i',j'}]=
\nu_i \nu_j \delta_{i,i'}\delta_{j,j'},
\end{eqnarray}
and other commutators vanish.
We impose the cocycle condition on 
the zero-mode operator $Q_{b}^{i,j}$ $(1\leq i<j \leq N+1)$ by
\begin{eqnarray}
~[Q_b^{i,j},Q_b^{i',j'}]=\delta_{j,N+1}\delta_{j',N+1}
\pi \sqrt{-1}~~~~~{\rm for}~(i,j) \neq (i',j').
\end{eqnarray}
We have the following (anti)commutation relations
\begin{eqnarray}
&&
\left[\exp\left(Q_b^{i,j}\right),\exp\left(Q_b^{i',j'}\right)
\right]=0
~~~
(1\leq i<j \leq N, 1\leq i'<j' \leq N),\\
&&\left\{\exp\left(Q_b^{i,N+1}\right),\exp\left(Q_b^{j,N+1}\right)
\right\}=0~~~
(1\leq i \neq j \leq N).
\end{eqnarray}
In what follows we use the standard normal ordering symbol $::$.
We set $b^{i,j}(z)$, $c^{i,j}(z)$, $b_\pm^{i,j}(z)$,
$a^j_\pm(z)$ and $
\left(\frac{\gamma_1}{\beta_1}\frac{\gamma_2}{\beta_2}
\cdots \frac{\gamma_r}{\beta_r}
~a^i \right)\left(z|\alpha \right)$ by
\begin{eqnarray}
&&a_\pm^{j}(z)=\pm (q-q^{-1})\sum_{\pm m>0}a_m^{j} 
z^{-m}\pm a_0^j {\rm log}q,\\
&&b^{i,j}(z)=
-\sum_{m \neq 0}\frac{b_m^{i,j}}{[m]}z^{-m}+Q_b^{i,j}+
b_0^{i,j}{\rm log}z,\\
&&b_\pm^{i,j}(z)=\pm (q-q^{-1})\sum_{\pm m>0}b_m^{i,j} 
z^{-m} \pm b_0^{i,j}{\rm log}q,\\
&&c^{i,j}(z)=
-\sum_{m \neq 0}\frac{c_m^{i,j}}{[m]}z^{-m}+Q_c^{i,j}+
c_0^{i,j}{\rm log}z,
\end{eqnarray}
\begin{eqnarray}
\left(\frac{\gamma_1}{\beta_1}\frac{\gamma_2}{\beta_2}
\cdots \frac{\gamma_r}{\beta_r}
~a^i \right)\left(z|\alpha \right)=
-\sum_{m \neq 0}\frac{[\gamma_1 m] \cdots [\gamma_r m]}
{[\beta_1 m] \cdots [\beta_r m]}
\frac{a^i_m}{[m]} q^{-\alpha |m|}z^{-m}
+\frac{\gamma_1 \cdots \gamma_r}{\beta_1 \cdots \beta_r}
(Q_a^i+a_0^i {\rm log}z).
\end{eqnarray}

\subsection{Bosonization of $U_q(\widehat{sl}(N|1))$}

We recall the bosonizations of the quantum superalgebra
$U_q(\widehat{sl}(N|1))$.

\begin{thm}
~~\cite{Kojima1}~~
A bosonization
of the quantum affine superalgebra
$U_q(\widehat{sl}(N|1))$ for an arbitrary level 
$k \in {\bf C}$
is given as follows.
For $1\leq i \leq N-1$ we set
\begin{eqnarray}
X_i^+(z)&=&\frac{1}{(q-q^{-1})z}
\sum_{j=1}^i(X_{i}^{+(j,1)}(z)-X_i^{+(j,2)}(z)),\\
X_N^+(z)&=& q^{N-2} \sum_{j=1}^N X_{N}^{+(j,0)}(z),
\label{boson2}
\\
X_i^-(z)&=&\frac{1}{(q-q^{-1})z}\left(
\sum_{j=1}^{i-1}(X_{i}^{-(j,1)}(z)-X_{i}^{-(j,2)}(z))
+(X_{i}^{-(i,1)}(z)-X_{i}^{-(i,2)}(z))
\right.
\nonumber\\
&&\left.-
\sum_{j=i+1}^{N-1}(X_{i}^{-(j,1)}(z)-X_{i}^{-(j,2)}(z))\right)
+q^{k+N-1}X_{i}^{-(N,0)}(z),
\label{boson3}
\\
X_N^-(z)&=&\frac{1}{(q-q^{-1})z}\sum_{j=1}^N q^{-N+j+1}
\left(-X_{N}^{-(j,1)}(z)+X_{N}^{-(j,2)}(z)\right).
\label{boson4}
\end{eqnarray}
\begin{eqnarray}
\Psi_i^\pm(q^{\pm \frac{k}{2}}z)&=&
\exp\left(
a_\pm^i(q^{\pm \frac{k+N-1}{2}}z)+
\sum_{l=1}^i(b_\pm^{l,i+1}(q^{\pm(l+k-1)}z)-b_\pm^{l,i}
(q^{\pm(l+k)}z))
\right.
\\
&&
\left.
+\sum_{l=i+1}^{N}(b_\pm^{i,l}(q^{\pm(k+l)}z)-
b_\pm^{i-1,l}(q^{\pm(k+l-1)}z))
+b_\pm^{i,N+1}(q^{\pm(k+N)}z)-
b_\pm^{i+1,N+1}(q^{\pm(k+N-1)}z)\right),\nonumber
\\
\Psi_N^\pm(q^{\pm \frac{k}{2}}z)&=&
\exp\left(a_\pm^N(q^{\pm \frac{k+N-1}{2}}z)-
\sum_{l=1}^{N-1}
(b_\pm^{l,N}(q^{\pm (k+l)}z)
+b_\pm^{l,N+1}(q^{\pm (k+l)}z))\right).
\label{boson6}
\end{eqnarray}
Here we have used
the auxiliary bosonic operators $X_{i}^{\pm (j,s)}(z)$ as follows.\\
For $1\leq i \leq N-1$ and $1\leq j \leq i$ we set
\begin{eqnarray}
X_{i}^{+ (j,1)}(z)&=&
:\exp\left((b+c)^{j,i}(q^{j-1}z)+b_+^{j,i+1}(q^{j-1}z)-
(b+c)^{j,i+1}(q^jz)\right.\nonumber\\
&&\left.+\sum_{l=1}^{j-1}
(b_+^{l,i+1}(q^{l-1}z)-b_+^{l,i}(q^lz))\right):,
\label{boson7}\\
X_{i}^{+ (j,2)}(z)&=&
:\exp\left((b+c)^{j,i}(q^{j-1}z)+b_-^{j,i+1}(q^{j-1}z)-
(b+c)^{j,i+1}(q^{j-2}z)\right.\nonumber\\
&&\left.+\sum_{l=1}^{j-1}
(b_+^{l,i+1}(q^{l-1}z)-b_+^{l,i}(q^lz))\right):.
\label{boson8}
\end{eqnarray}
For $1\leq j \leq N$ we set
\begin{eqnarray}
X_{N}^{+ (j,0)}(z)&=&:\exp\left(
(b+c)^{j,N}(q^{j-1}z)
+b^{j,N+1}(q^{j-1}z)
-\sum_{l=1}^{j-1}(b_+^{l,N+1}(q^lz)+b_+^{l,N}(q^lz))
\right):.\label{boson9}
\end{eqnarray}
For $1\leq i \leq N-1$ and $1\leq j \leq i-1$ we set
\begin{eqnarray}
X_{i}^{- (j,1)}(z)&=&
:\exp\left(
a_-^i(q^{-\frac{k+N-1}{2}}z)
+(b+c)^{j,i+1}(q^{-k-j}z)-b_-^{j,i}(q^{-k-j}z)
-(b+c)^{j,i}(q^{-k-j+1}z)\right.\nonumber\\
&& 
+\sum_{l=j+1}^i 
(b_-^{l,i+1}(q^{-k-l+1}z)-b_-^{l,i}(q^{-k-l}z))
+\sum_{l=i+1}^N
(b_-^{i,l}(q^{-k-l}z)-b_-^{i+1,l}(q^{-k-l+1}z))\nonumber\\
&&
\left.+b_-^{i,N+1}(q^{-k-N}z)-b_-^{i+1,N+1}(q^{-k-N+1}z)
\right):,
\label{boson10}
\\
X_{i}^{- (j,2)}(z)&=&
:\exp\left(a_-^i(q^{-\frac{k+N-1}{2}}z)
+(b+c)^{j,i+1}(q^{-k-j}z)-b_+^{j,i}(q^{-k-j}z)
-(b+c)^{j,i}(q^{-k-j-1}z)\right.
\nonumber\\
&& 
+\sum_{l=j+1}^i 
(b_-^{l,i+1}(q^{-k-l+1}z)-b_-^{l,i}(q^{-k-l}z))
+\sum_{l=i+1}^N
(b_-^{i,l}(q^{-k-l}z)-b_-^{i+1,l}(q^{-k-l+1}z))\nonumber\\
&&
\left.+b_-^{i,N+1}(q^{-k-N}z)-b_-^{i+1,N+1}(q^{-k-N+1}z)\right):.
\label{boson11}
\end{eqnarray}
For $1\leq i \leq N-1$ we set
\begin{eqnarray}
X_{i}^{- (i,1)}(z)&=&:\exp\left(a_-^i(q^{-\frac{k+N-1}{2}}z)
+(b+c)^{i,i+1}(q^{-k-i}z)+\sum_{l=i+1}^N(b_-^{i,l}(q^{-k-l}z)
-b_-^{i+1,l}(q^{-k-l+1}z))\right.
\nonumber\\
&&
\left.
+b_-^{i,N+1}(q^{-k-N}z)-b_-^{i+1,N+1}(q^{-k-N+1}z)\right):,
\label{boson12}
\\
X_{i}^{- (i,2)}(z)&=&:
\exp\left(a_+^i(q^{\frac{k+N-1}{2}}z)
+(b+c)^{i,i+1}(q^{k+i}z)+\sum_{l=i+1}^N(b_+^{i,l}(q^{k+l}z)
-b_+^{i+1,l}(q^{k+l-1}z))\right.
\nonumber\\
&&
\left.
+b_+^{i,N+1}(q^{k+N}z)-b_+^{i+1,N+1}(q^{k+N-1}z)\right):.
\label{boson13}
\end{eqnarray}
For $1\leq i \leq N-1$ and $i+1 \leq j \leq N-1$ we set
\begin{eqnarray}
X_{i}^{- (j,1)}(z)&=&
:\exp\left(a_+^i(q^{\frac{k+N-1}{2}}z)
+(b+c)^{i,j+1}(q^{k+j}z)
+b_+^{i+1,j+1}(q^{k+j}z)-(b+c)^{i+1,j+1}(q^{k+j+1}z)\right.
\nonumber\\
&&
\left.
+\sum_{l=j+1}^N
(b_+^{i,l}(q^{k+l}z)-b_+^{i+1,l}(q^{k+l-1}z))+b_+^{i,N+1}(q^{k+N}z)-b_+^{i+1,N+1}(q^{k+N-1}z)\right):,
\label{boson14}
\\
X_{i}^{- (j,2)}(z)&=&
:\exp\left(a_+^i(q^{\frac{k+N-1}{2}}z)
+(b+c)^{i,j+1}(q^{k+j}z)
+b_-^{i+1,j+1}(q^{k+j}z)-(b+c)^{i+1,j+1}(q^{k+j-1}z)\right.
\nonumber\\
&&
\left.+\sum_{l=j+1}^N
(b_+^{i,l}(q^{k+l}z)-b_+^{i+1,l}(q^{k+l-1}z))
+b_+^{i,N+1}(q^{k+N}z)-b_+^{i+1,N+1}(q^{k+N-1}z)\right):.
\label{boson15}
\end{eqnarray}
For $1\leq i \leq N-1$ we set
\begin{eqnarray}
X_{i}^{- (N,0)}(z)&=&:\exp\left(a_+^i(q^{\frac{k+N-1}{2}}z)
-b^{i,N+1}(q^{k+N-1}z)-b_+^{i+1,N+1}(q^{k+N-1}z)+b^{i+1,N+1}(q^{k+N}z)
\right):.\nonumber\\
\label{boson16}
\end{eqnarray}
For $1\leq j \leq N-1$ we set
\begin{eqnarray}
X_{N}^{- (j,1)}(z)&=&:
\exp\left(a_-^N(q^{-\frac{k+N-1}{2}}z)
-b_-^{j,N}(q^{-k-j}z)-(b+c)^{j,N}(q^{-k-j+1}z)\right.\\
&&
\left.
-b_-^{j,N+1}(q^{-k-j}z)-b^{j,N+1}(q^{-k-j+1}z)
-\sum_{l=j+1}^{N-1}
(b_-^{l,N}(q^{-k-l}z)+b_-^{l,N+1}(q^{-k-l}z))\right):,
\nonumber
\\
X_{N}^{- (j,2)}(z)&=&:\exp\left(
a_-^N(q^{-\frac{k+N-1}{2}}z)
-b_+^{j,N}(q^{-k-j}z)-(b+c)^{j,N}(q^{-k-j-1}z)\right.
\\
&&\left.
-b_+^{j,N+1}(q^{-k-j}z)-b^{j,N+1}(q^{-k-j-1}z)
-\sum_{l=j+1}^{N-1}
(b_-^{l,N}(q^{-k-l}z)+b_-^{l,N+1}(q^{-k-l}z))\right):,
\nonumber
\\
X_{N}^{- (N,1)}(z)&=&
:\exp\left(
a_-^N(q^{-\frac{k+N-1}{2}}z)-b^{N,N+1}(q^{-k-N+1}z)\right):,
\label{boson19}
\\
X_{N}^{- (N,2)}(z)&=&:
\exp\left(
a_+^N(q^{\frac{k+N-1}{2}}z)-b^{N,N+1}(q^{k+N-1}z)\right):.
\label{boson20}
\end{eqnarray}
The ${\bf Z}_2$-grading is :
$|X_{N}^{\pm (j,s)}(z)|=1$ and zero otherwise.
\end{thm}
Very explicitly, we have
\begin{eqnarray}
h_{i,m}&=&
q^{-\frac{N-1}{2}|m|}a_m^i+\sum_{l=1}^i
(q^{-(\frac{k}{2}+l-1)|m|}b_m^{l,i+1}-
q^{-(\frac{k}{2}+l)|m|}b_m^{l,i})\nonumber\\
&&+\sum_{l=i+1}^N(q^{-(\frac{k}{2}+l)|m|}b_m^{i,l}-
q^{-(\frac{k}{2}+l-1)|m|}b_m^{i+1,l})\nonumber\\
&&+q^{-(\frac{k}{2}+N)|m|}b_m^{i,N+1}-
q^{-(\frac{k}{2}+N-1)|m|}b_m^{i+1,N+1}~~~~~(1\leq i \leq N-1),\\
h_{N,m}&=&q^{-\frac{N-1}{2}|m|}a_m^N-
\sum_{l=1}^{N-1}(q^{-(\frac{k}{2}+l)|m|}b_m^{l,N}+
q^{-(\frac{k}{2}+l)|m|}b_m^{l,N+1}).
\end{eqnarray}

\subsection{Fock-Wakimoto module}

We introduce the vacuum state $|0\rangle$ of 
the boson Fock space by
\begin{eqnarray}
a_m^i|0\rangle=b_m^{i,j}|0\rangle
=c_m^{i,j}|0 \rangle=0~~(m \geq 0).
\end{eqnarray}
For $p_a^i \in {\bf C}$ $(1\leq i \leq N)$,
$p_b^{i,j} \in {\bf C}$ $(1\leq i<j \leq N+1)$,
$p_c^{i,j} \in {\bf C}$ $(1\leq i<j \leq N)$,
we set
\begin{eqnarray}
&&|p_a, p_b, p_c \rangle
\\
&=&\exp\left(
\sum_{i,j=1}^N
\frac{{\rm Min}(i,j)(N-1-{\rm Max}(i,j))}{(N-1)(k+N-1)}
p_a^i Q_a^j
-\sum_{1\leq i<j \leq N+1}p_b^{i,j}Q_b^{i,j}
+\sum_{1\leq i<j \leq N}p_c^{i,j}Q_c^{i,j}
\right)|0\rangle.\nonumber
\end{eqnarray}
It satisfies
\begin{eqnarray}
a_0^i|p_a,p_b,p_c\rangle=p_a^i |p_a,p_b,p_c\rangle,~
b_0^{i,j}|p_a,p_b,p_c\rangle=p_b^{i,j} |p_a,p_b,p_c\rangle,~
c_0^{i,j}|p_a,p_b,p_c\rangle=p_c^{i,j} |p_a,p_b,p_c\rangle.
\end{eqnarray}
The boson Fock space $F(p_a,p_b,p_c)$
is generated by
the bosons $a_m^i, b_m^{i,j}, c_m^{i,j}$
on the vector $|p_a,p_b,p_c\rangle$.
We set the space $F(p_a)$ by
\begin{eqnarray}
F(p_a)=
\bigoplus
_{
p_b^{i,j}=-p_c^{i,j} \in {\bf Z}~(1\leq i<j \leq N)
\atop{
p_b^{i,N+1} \in {\bf Z}~(1\leq i \leq N)
}}F(p_a,p_b,p_c).
\end{eqnarray}
We impose the restriction
$p_b^{i,j}=-p_c^{i,j} \in {\bf Z}$ $(1\leq i<j \leq N)$,
because 
the $X_{i,m}^\pm$ change 
$Q_b^{i,j}+Q_c^{i,j}$.
The $F(p_a)$ is $U_q(\widehat{sl}(N|1))$-module.
We set the vector $|\lambda\rangle=|p_a, 0, 0\rangle$
upon the specialization $p_b^{i,j}=0~(1\leq i<j \leq N+1)$ 
and $p_c^{i,j}=0~(1\leq i<j \leq N)$.

\begin{prop}~~~
\label{prop:space1}
The $|\lambda\rangle=|p_a, 0, 0 \rangle$
is the highest weight vector of the highest weight
whose classical part is
$\bar{\lambda}=\sum_{j=1}^N
p_a^j \bar{\Lambda}_j$.
\begin{eqnarray}
&&h_{i,m}|\lambda\rangle=0,~~~
X_{i,m}^\pm |\lambda \rangle=0,~~~(m>0),\\
&&X_{i,0}^+ |\lambda \rangle=0,~~~
h_{i,0}|\lambda \rangle=p_a^i |\lambda \rangle.
\end{eqnarray}
Using the highest weight vector $|\lambda \rangle$,
we have the highest weight module $V(\lambda)$ of
$U_q(\widehat{sl}(N|1))$.
\begin{eqnarray}
V(\lambda)\subset F(p_a).
\end{eqnarray}
\end{prop}
The module $F(p_a)$
is not irreducible.
We recall the non-quantum algebra 
$\widehat{sl}(2)$ case \cite{Bernard-Felder}. 
The irreducible highest weight module $L(\lambda)$
for the affine algebra $\widehat{sl}(2)$
was constructed from the Fock-Wakimoto module
on the boson Fock space \cite{Wakimoto}
by the Felder complex. 
We recall the quantum algebra $U_q(\widehat{sl}(2))$ 
case \cite{Matsuo, Shiraishi, Konno}.
The irreducible highest weight
module $L(\lambda)$ for 
$U_q(\widehat{sl}(2))$
was constructed from the similar space as $F(p_a)$
by two steps; the first step is
the construction of the Fock-Wakimoto module by the $\xi$-$\eta$ system, and
the second step is the resolution
by the Felder complex \cite{Konno}.
The submodule of the quantum algebra $U_q(\widehat{sl}(2))$,
induced by the $\xi$-$\eta$ system,
plays the same role as the Fock-Wakimoto module
of the non-quantum algebra $\widehat{sl}(2)$.
We call this submodule
induced by the $\xi$-$\eta$ system
the "Fock-Wakimoto module".
In this paper we study the $\xi$-$\eta$ system
and propose the Fock-Wakimoto module 
for $U_q(\widehat{sl}(N|1))$.

\begin{dfn}~~
We introduce the operators $\xi_m^{i,j}$ and
$\eta_m^{i,j}$ $(1\leq i<j \leq N, m \in {\bf Z})$ by
\begin{eqnarray}
\eta^{i,j}(z)=\sum_{m \in {\bf Z}}\eta_{m}^{i,j}
z^{-m-1}=:e^{c^{i,j}(z)}:,~~
\xi^{i,j}(z)=\sum_{m \in {\bf Z}}\xi_{m}^{i,j}z^{-m}=:e^{-c^{i,j}(z)}:.
\end{eqnarray}
\end{dfn}

The Fourier components $\eta_m^{i,j}
=\oint \frac{dz}{2\pi \sqrt{-1}}z^m \eta^{i,j}(z)$,
$\xi_m^{i,j}=
\oint \frac{dz}{2\pi \sqrt{-1}}z^{m-1}\xi^{i,j}(z)$ $(m \in {\bf Z})$
are well defined on the space ${F}(p_a)$.
They satisfy the anti-commutation relations.
\begin{eqnarray}
\{\eta_m^{i,j},\xi_n^{i,j}\}=\delta_{m+n,0},~
\{\eta_m^{i,j},\eta_n^{i,j}\}=\{\xi_m^{i,j},\xi_n^{i,j}\}=0~~~
(1\leq i<j \leq N).
\end{eqnarray}
They commute with each other
\begin{eqnarray}
[\eta_m^{i,j},\xi_n^{i',j'}]=[\eta_m^{i,j},\eta_n^{i',j'}]
=[\xi_m^{i,j},\xi_n^{i',j'}]=0~~~(i,j)\neq (i',j').
\end{eqnarray}
We focus our attention on the operators 
$\eta_0^{i,j}, \xi_0^{i,j}$ satisfying
$(\eta_0^{i,j})^2=0$, $(\xi_0^{i,j})^2=0$.
They satisfy
\begin{eqnarray}
{\rm Im} (\eta_0^{i,j})={\rm Ker} (\eta_0^{i,j}),~~~~~
{\rm Im} (\xi_0^{i,j})={\rm Ker} (\xi_0^{i,j}).
\end{eqnarray}
The products $\eta_0^{i,j} \xi_0^{i,j}$ and 
$\xi_0^{i,j} \eta_0^{i,j}$ are 
the projection operators, which satisfy
\begin{eqnarray}
\eta_0^{i,j}\xi_0^{i,j}+\xi_0^{i,j}\eta_0^{i,j}=1,
\end{eqnarray}
and
\begin{eqnarray}
(\eta_0^{i,j}\xi_0^{i,j})^2=\eta_0^{i,j}\xi_0^{i,j},
~(\xi_0^{i,j}\eta_0^{i,j})^2=\xi_0^{i,j}\eta_0^{i,j},
~(\xi_0^{i,j}\eta_0^{i,j})(\eta_0^{i,j}\xi_0^{i,j})=0,
~(\eta_0^{i,j}\xi_0^{i,j})(\xi_0^{i,j}\eta_0^{i,j})=0.
\end{eqnarray}
Hence we have a direct sum decomposition.
\begin{eqnarray}
F(p_a)=
\eta_0^{i,j}\xi_0^{i,j}F(p_a) \oplus 
\xi_0^{i,j}\eta_0^{i,j}F(p_a),
\end{eqnarray}
and
\begin{eqnarray}
{\rm Ker} (\eta_0^{i,j})=
\eta_0^{i,j}\xi_0^{i,j}F(p_a),~~~~~
{\rm Coker} (\eta_0^{i,j})=
\xi_0^{i,j}\eta_0^{i,j}F(p_a)
=F_\Lambda / 
(\eta_0^{i,j}\xi_0^{i,j})F_\Lambda.
\end{eqnarray}
We set 
\begin{eqnarray}
\eta_0=\prod_{1\leq i < j \leq N}\eta_0^{i,j},~~~~~
\xi_0=\prod_{1\leq i < j \leq N}\xi_0^{i,j}.
\end{eqnarray}

\begin{dfn}~~~
We introduce the subspace ${\cal F}(p_a)$ by 
\begin{eqnarray}
{\cal F}(p_a)
=\eta_0 \xi_0
F(p_a).
\end{eqnarray}
\end{dfn}

The operators $\eta_0^{i,j}$, $\xi_0^{i,j}$ commute with
the operators $X_{i',j'}^\pm(z)$, $\Psi_{i'}^\pm(z)$ up to sign $\pm$.
When we set 
the operators $\widetilde{X}_{i}^\pm(z)$, $\widetilde{\Psi}_i^\pm(z)$ by
the conditions
$\widetilde{X}_{i}^\pm(z) \eta_0^{i',j'}=
\eta_0^{i',j'}{X}_i^\pm(z)$, $\widetilde{\Psi}_i^\pm(z)\eta_0^{i',j'}=
\eta_0^{i',j'} \Psi_i^\pm(z)$,
the bosonic operators $\widetilde{X}_{i}^\pm(z)$, $\widetilde{\Psi}_i^\pm(z)$ 
give a bosonization of $U_q(\widehat{sl}(N|1))$ again.

\begin{prop}~\cite{Kojima2}~~
The subspace ${\cal F}(p_a)$ is the $U_q(\widehat{sl}(N|1))$ module.
\end{prop}
We call the submodule ${\cal F}(p_a)$ the Fock-Wakimoto module.
It is expected that we have
the irreducible highest weight module $L(\lambda)$
with the highest weight $\lambda$,
whose classical part 
$\bar{\lambda}=\sum_{j=1}^N p_a^i \bar{\Lambda}_i$,
by the Felder complex.
The construction of the Felder complex
is open problem even for non-superalgebra 
$U_q(\widehat{sl}(3))$. We would like to report
the Felder complex of $U_q(\widehat{sl}(N))$ and 
$U_q(\widehat{sl}(N|1))$ in the future publications.

\section{Screening current}

In this section we introduce the screening operators 
$S_i$ $(i=1,2,\cdots,N)$,
which commute with $U_q(\widehat{sl}(N|1))$ for
an arbitrary level $k \neq -N+1$.
We need the screening operators
to construct the vertex operators.

\subsection{Screening current}

We set the $q$-difference operators with
a parameter $\alpha$ by
\begin{eqnarray}
(_\alpha \partial_z f)(z)=\frac{
f(q^{\alpha} z)-f(q^{-\alpha} z)}{(q-q^{-1})z}.
\end{eqnarray}
The Jackson integral with parameter $p \in {\bf C}$ $(|p|<1)$
and $s \in {\bf C}^*$ is defined by
\begin{eqnarray}
\int_0^{s \infty}f(z)d_pz=s(1-p)\sum_{m \in {\bf Z}}f(sp^m)p^m.
\end{eqnarray}
The Jackson integral satisfies 
\begin{eqnarray}
\int_0^{s \infty}
(_\alpha \partial_z f)(z)
d_{p}z=0~~~(p=q^{2\alpha}).
\end{eqnarray}
For $r \in {\bf C}$ $({\rm Re}(r)>0)$ we introduce
the Jacobi elliptic theta function
\begin{eqnarray}
~[u]_r=
q^{\frac{u^2}{r}-u}\frac{\Theta_{q^{2r}}(q^{2u})}{
(q^{2r};q^{2r})_\infty^3},
\end{eqnarray}
where we have used
\begin{eqnarray}
\Theta_p(z)=(z;p)_\infty 
(pz^{-1};p)_\infty (p;p)_\infty,~~~
(z;p)_\infty
=\prod_{m=0}^\infty (1-p^{m}z).
\end{eqnarray}
The Jacobi elliptic theta function
satisfies the quasi-periodicity property
\begin{eqnarray}
~[u+r]_r=-[u]_r,~~
~[u+r\tau]_r=-e^{-\pi i \tau-\frac{2\pi i}{r}u}[u]_r,
\end{eqnarray}
where $\tau$ such that
${\rm Im}(\tau)>0$ is given by $q^{2r}=e^{-\frac{2\pi i}{\tau}}$.

\begin{dfn}~~~We introduce 
the bosonic operators $S_i(z)$ $(i=1,2,\cdots, N)$ 
that we call the screening current as follows.
\begin{eqnarray}
S_i(z)&=&\frac{1}{(q-q^{-1})z}\sum_{j=i+1}^N
(S_{i}^{ (j,1)}(z)-S_{i}^{ (j,2)}(z))+q S_{i}^{ (N+1,0)}(z)
~~~~~
(1\leq i \leq N-1),\\
S_N(z)&=&-q^{-1}S_{N}^{ (N+1,0)}(z).
\end{eqnarray}
For $2i+1 \leq j \leq 2N+1$ we have set
\begin{eqnarray}
S_{i}^{ (j,s)}(z)
=:\exp\left(-\left(\frac{1}{k+N-1}a^i\right)
\left(z\left|\frac{k+N-1}{2}\right.\right)\right)
\widetilde{S}_{i}^{ (j,s)}(z):.
\end{eqnarray}
Here, for $1\leq i \leq N-1$ and $i+1 \leq j \leq N$,
we have set
\begin{eqnarray}
\widetilde{S}_{i}^{ (j,1)}(z)&=&
:\exp\left(-b_-^{i,j}(q^{N-1-j}z)-(b+c)^{i,j}(q^{N-j}z)
+(b+c)^{i+1,j}(q^{N-1-j}z)\right.\nonumber\\
&&\left.+\sum_{l=j+1}^N
(b_-^{i+1,l}(q^{N-l}z)-b_-^{i,l}(q^{N-l-1}z))
+b_-^{i+1,N+1}(z)-b_-^{i,N+1}(q^{-1}z)\right):,\\
\widetilde{S}_{i}^{ (j,2)}(z)&=&
:\exp\left(-b_+^{i,j}(q^{N-1-j}z)-(b+c)^{i,j}(q^{N-j-2}z)
+(b+c)^{i+1,j}(q^{N-1-j}z)\right.\nonumber\\
&&\left.+
\sum_{l=j+1}^N
(b_-^{i+1,l}(q^{N-l}z)-b_-^{i,l}(q^{N-l-1}z))
+b_-^{i+1,N+1}(z)-b_-^{i,N+1}(q^{-1}z)\right):.
\end{eqnarray}
For $1\leq i \leq N-1$ we have set
\begin{eqnarray}
\widetilde{S}_{i}^{ (N+1,0)}(z)&=&:\exp\left(b^{i,N+1}(z)
+b_+^{i+1,N+1}(z)-b^{i+1,N+1}(qz)\right):,
\\
\widetilde{S}_{N}^{ (N+1,0)}(z)&=&:\exp\left(b^{N,N+1}(z)\right):.
\end{eqnarray}
The ${\bf Z}_2$-grading of the screening currents are :
$|S_{N}^{(N+1,0)}(z)|=1$ and zero otherwise.
\end{dfn}

\begin{thm}
\label{thm:screening}
~~~The screening currents $S_i(z)$ $(i=1,2,\cdots, N)$ commute 
(or anti-commute) with
$U_q(\widehat{sl}(N|1))$ modulo total difference.
\begin{eqnarray}
~[h_{i,m},S_j(z)]&=&0,\label{thm:screening1}\\
~[X_i^+(z_1),S_j(z_2)]&=&0,\label{thm:screening2}\\
~[X_i^-(z_1),S_j(z_2)]&=&\frac{\delta_{i,j}}{(q-q^{-1})z_1^2}
\left(_{k+N-1}\partial_{z}\delta \right)(z_2/z_1)\nonumber\\
&\times&
:\exp\left(-\left(\frac{1}{k+N-1}~a^j\right)\left(z_1\left|
-\frac{k+N-1}{2}\right.\right)\right):.
\label{thm:screening3}
\end{eqnarray}
The screening currents $S_i(z)$ $(i=1,2,\cdots,N)$ 
satisfy
\begin{eqnarray}
\left[u_1-u_2+\frac{A_{i,j}}{2}\right]_{k+N-1}S_i(z_1)S_j(z_2)&=&
\left[u_2-u_1+\frac{A_{i,j}}{2}\right]_{k+N-1}S_j(z_1)S_i(z_2).
\label{thm:screening4}
\end{eqnarray}
The symbol $[u]_{k+N-1}$ represents 
the Jacobi elliptic theta function.
Here we have used $z_j=q^{2u_j}$.
\end{thm}

\begin{dfn}~~~
We introduce the screening operators 
$Q_i$ $(i=1,2,\cdots,N)$ 
by the Jackson integral.
\begin{eqnarray}
Q_i=\int_0^{s \infty}S_i(z)d_pz,~~~(p=q^{2(k+N-1)}).
\end{eqnarray}
The screening operators $Q_i$ are convergent 
on the Fock space.
\end{dfn}

\begin{cor}~~~
The screening operators $Q_i$ $(i=1,2,\cdots,N)$
commute with the quantum superalgebra 
$U_q(\widehat{sl}(N|1))$.
\end{cor}

\begin{prop}~~~
The screening operators $Q_i$ $(i=1,2,\cdots,N)$
commute with the projection operator 
$\eta_0 \xi_0$ of the $\xi$-$\eta$ system.
Hence the screening operators $Q_i$
act on the Fock-Wakimoto module ${\cal F}(p_a)$.
\end{prop}

\subsection{Proof}

Here we give proof of theorem \ref{thm:screening}.
Direct calculations of the normal orderings
show theorem \ref{thm:screening}.

~\\
$\bullet$
~Proof of (\ref{thm:screening3}) for $1\leq i=j \leq N$.
\\
First we show (\ref{thm:screening3}) 
for $1\leq i=j \leq N-1$.
The commutators vanish,
$[X_i^{- (l,s)}(z_1),S_i^{(m,t)}(z_2)]=0$,
for the following condition.
\begin{eqnarray}
~[(l,s),(m,t)]\neq 
\left\{\begin{array}{cc}
~[(i,1),(i+1,1)],
~[(i,2),(i+1,2)]
&~~(1\leq i \leq N-1)\\
~[(l,1),(l+1,2)],
~[(l,2),(l+1,1)]
&~(1\leq i \leq N-1, i+1\leq l\leq N-1)\\
~[(N,0),(N+1,0)] & ~(1\leq i \leq N-1)
\end{array}\right.
.\label{proof1}
\end{eqnarray}
Hence we have the following relations
for $1\leq i \leq N-1$.
\begin{eqnarray}
~[X_i^-(z_1),S_i(z_2)]&=&q^{k+N}[X_i^{- (N,0)}(z_1),S_i^{(N+1,0)}(z_2)]
\nonumber\\
&&+\frac{1}{(q-q^{-1})^2z_1z_2}
\left\{[X_i^{- (i,1)}(z_1),S_i^{(i+1,1)}(z_2)]
+[X_i^{- (i,2)}(z_1),S_i^{(i+1,2)}(z_2)]\right.\nonumber\\
&&\left.+\sum_{l=i+1}^{N-1}
\left([X_i^{- (l,1)}(z_1),S_i^{(l+1,2)}(z_2)]
+[X_i^{- (l,2)}(z_1),S_i^{(l+1,1)}(z_2)]\right)
\right\}.
\end{eqnarray}
Using the relations (\ref{A1-1}), (\ref{A1-2}), 
(\ref{A1-3}), (\ref{A1-4}), 
(\ref{A1-5}), (\ref{A1-6}) in appendix \ref{appendixA},
we have
\begin{eqnarray}
&&[X_i^-(z_1),S_i(z_2)]
(q-q^{-1})z_1z_2
\nonumber\\
&=&\delta\left(\frac{q^{-N-k-1}z_2}{z_1}\right)
(:X_i^{- (N,0)}(z_1)S_i^{(N+1,0)}(z_2):
-:X_i^{- (N-1,1)}(z_1)S_i^{(N,2)}(z_2):)\nonumber\\
&+&\delta\left(\frac{q^{N-3-k-2i}z_2}{z_1}\right)
(:X_i^{- (i+1,2)}(z_1)S_i^{(i+2,1)}(z_2):
-:X_i^{- (i,2)}(z_1)S_i^{(i+1,2)}(z_2):)\\
&+&\sum_{l=i+1}^{N-2}
\delta\left(\frac{q^{N-3-k-2l}z_2}{z_1}\right)
(:X_i^{- (l+1,2)}(z_1)S_i^{(l+2,1)}(z_2):
-:X_i^{- (l+1,1)}(z_1)S_i^{(l+2,2)}(z_2):)\nonumber\\
&+&\delta\left(\frac{q^{N+k-1}z_2}{z_1}\right)
:X_i^{- (i,1)}(z_1)S_i^{(i+1,1)}(z_2):
-
\delta\left(\frac{q^{-N-k+1}z_2}{z_1}\right)
:X_i^{- (N,0)}(z_1)S_i^{(N+1,0)}(z_2):.\nonumber
\end{eqnarray}
Using specializations (\ref{A2-1}), (\ref{A2-3}), (\ref{A2-4}), (\ref{A2-5})
in appendix \ref{appendixA},
we conclude (\ref{thm:screening3}) for $1\leq i=j \leq N-1$.
Next we show (\ref{thm:screening3}) 
for $i=j=N$.
The commutators vanish,
$[X_N^{- (l,s)}(z_1),S_N^{(m,t)}(z_2)]=0$,
for $[(l,s),(m,t)]
\neq [(N,1), (N+1,0)], [(N,2), (N+1,0)]$.
Hence, using the relations (\ref{A1-7}) and (\ref{A1-8}) 
in appendix \ref{appendixA}, we have
\begin{eqnarray}
&&~[X_N^-(z_1),S_N(z_2)](q-q^{-1})z_1z_2\\
&=&
\delta\left(\frac{q^{N+k-1}z_2}{z_1}\right)
:X_N^{- (N,1)}(z_1)S_N^{(N+1,0)}(z_2):-
\delta\left(\frac{q^{N+k-1}z_2}{z_1}\right)
:X_N^{- (N,2)}(z_1)S_N^{(N+1,0)}(z_2):.\nonumber
\end{eqnarray}
Using the relation (\ref{A2-2}) 
in appendix \ref{appendixA},
we have (\ref{thm:screening3}) for $i=N$.
Now we have shown (\ref{thm:screening3}) for $1\leq i=j \leq N$.
\\
$\bullet$~Proof of (\ref{thm:screening3}) for
$1\leq i \neq j \leq N$.
\\
First we show (\ref{thm:screening3}) for $i+1<j$.
In this case the commutators vanish,
$[X_i^{- (l,s)}(z_1),S_j^{(m,t)}(z_2)]=0$, for every $[(l,s),(m,t)]$.
Hence we conclude $[X_i^-(z_1),S_j(z_2)]=0$.
Next we show (\ref{thm:screening3}) for $j=i+1$.
The commutators vanish,
$[X_i^{- (l,s)}(z_1),S_j^{(m,t)}(z_2)]=0$, for 
the following condition.
\begin{eqnarray}
[(l,s),(m,t)]\neq [(l,1), (l+1,2)], [(l,2), (l+1,1)]
~~(1\leq i \leq N-2, i+1\leq l \leq N-1).
\end{eqnarray}
Hence, using the relations (\ref{A1-9}) and (\ref{A1-10}), 
we have the following relation
for $1\leq i \leq N-1$.
\begin{eqnarray}
&&[X_i^-(z_1),S_{i+1}(z_2)](q-q^{-1})z_1z_2\\
&=&
\sum_{l=i+1}^{N-1}\delta\left(\frac{q^{N-k-2l-2}z_2}{z_1}\right)
\left(-:X_i^{- (l,1)}(z_1)S_{i+1}^{(l+1,2)}(z_2):+
:X_i^{- (l,2)}(z_1)S_{i+1}^{(l+1,1)}(z_2):\right).
\nonumber
\end{eqnarray}
Using the specialization
(\ref{A2-6}), we conclude (\ref{thm:screening3}) for $j=i+1$.
Next we show (\ref{thm:screening3}) for $1\leq j<i \leq N-1$.
The commutators vanish,
$[X_i^{- (l,s)}(z_1),S_j^{(m,t)}(z_2)]=0$, for 
the following condition.
\begin{eqnarray}
~[(l,s),(m,t)]\neq
\left\{
\begin{array}{cc}
\begin{array}{c}
~[(j,1), (i+1,1)], ~[(j+1,1), (i,1)]\\
~[(j,2), (i+1,1)], ~[(j+1,1), (i,2)]\\
~[(j,2), (i,1)], ~[(j,1), (i,2)]
\end{array}
&(1\leq j<i \leq N-1)\\
\end{array}
\right..
\end{eqnarray}
Hence, using the relations
(\ref{A1-11}), (\ref{A1-12}), (\ref{A1-13}), (\ref{A1-14}), 
(\ref{A1-15}), (\ref{A1-16}),
we have the following relation for
$1\leq j<i \leq N-1$.
\begin{eqnarray}
&&[X_i^-(z_1),S_j(z_2)](q-q^{-1})z_1z_2\nonumber\\
&=&\delta\left(\frac{q^{N+k-i+j-1}z_2}{z_1}\right)
\left(:X_i^{- (j,1)}(z_1)S_j^{(i+1,1)}(z_2):
-:X_i^{- (j+1,1)}(z_1)S_j^{(i,1)}(z_2):\right.\nonumber\\
&&
-:X_i^{- (j,2)}(z_1)S_j^{(i+1,1)}(z_2):
+:X_i^{- (j+1,1)}(z_1)S_j^{(i,2)}(z_2):\nonumber\\
&&\left.-q^{-1}:X_i^{- (j,2)}(z_1)S_j^{(i,1)}(z_2):
+q^{-1}:X_i^{- (j,1)}(z_1)S_j^{(i,2)}(z_2):
\right).
\end{eqnarray}
Using the specializations 
(\ref{A2-7}), (\ref{A2-8}), (\ref{A2-9}), we have
$[X_i^-(z_1),S_j(z_2)]=0$ for $1\leq j <i \leq N-1$.
Next we show (\ref{thm:screening3}) for $1\leq j \leq N-1$
and $i=N$.
The commutators vanish,
$[X_i^{- (l,s)}(z_1),S_j^{(m,t)}(z_2)]=0$, for 
the following condition.
\begin{eqnarray}
~[(l,s),(m,t)]\neq
\left\{
\begin{array}{cc}
\begin{array}{c}
~[(j,1), (N,2)], ~[(j,2), (N,1)]\\
~[(j,1), (N+1,0)], ~[(j+1,1), (N,1)]\\
~[(j,2), (N+1,0)], ~[(j+1,1), (N,2)]
\end{array}
~~(1\leq j \leq N-1)
\end{array}
\right..
\end{eqnarray}
Hence, using the relations (\ref{A1-12}), (\ref{A1-14}), 
(\ref{A1-15}), (\ref{A1-16}), (\ref{A1-17}), 
(\ref{A1-18}) we have
the following relation for
$1\leq j \leq N-1$.
\begin{eqnarray}
&&[X_N^-(z_1),S_j(z_2)](q-q^{-1})q^{N-j-2}z_1z_2\nonumber\\
&=&\delta\left(\frac{q^{k+j-1}z_2}{z_1}\right)
\left(q^{-2}:X_N^{- (j,1)}(z_1)S_j^{(N,2)}(z_2):
-q^{-2}:X_N^{- (j,2)}(z_1)S_j^{(N,1)}(z_2):\right.\nonumber\\
&&
-:X_N^{-(j,1)}(z_1)S_j^{(N+1,0)}(z_2):+
:X_N^{- (j+1,1)}(z_1)S_j^{(N,1)}(z_2):\nonumber\\
&&\left.+
-:X_N^{-(j,2)}(z_1)S_j^{(N+1,0)}(z_2):+
:X_N^{- (j+1,1)}(z_1)S_j^{(N,2)}(z_2):
\right).
\end{eqnarray}
Using the specializations 
(\ref{A2-9}), (\ref{A2-10}), (\ref{A2-11}), we have
$[X_N^-(z_1),S_j(z_2)]=0$ for $1\leq j \leq N-1$.
Now we have shown $[X_i^-(z_1),S_j(z_2)]=0$ for $1\leq i \neq j \leq N$.
\\
$\bullet$~Proof of (\ref{thm:screening2}) 
for $1\leq i=j \leq N$.
\\
The commutators vanish,
$[X_i^{+ (l,s)}(z_1),S_i^{(m,t)}(z_2)]=0$,
for the condition
\begin{eqnarray}
[(l,s),(m,t)]
\neq [(i,1), (i+1,2)], [(i,2), (i+1,1)]~~~
(1\leq i \leq N-1).
\end{eqnarray}
We have $[X_N^+(z_1),S_N(z_2)]=0$.
For $1\leq i \leq N-1$, using the relations 
(\ref{A1-19}) and (\ref{A1-20}), we have
\begin{eqnarray}
&&[X_i^+(z_1), S_i(z_2)](q-q^{-1})z_1z_2\nonumber\\
&=&\delta\left(\frac{q^{N-2i-1}z_2}{z_1}\right)
\left(-:X_i^{+ (i,1)}(z_1)S_i^{(i+1,2)}(z_2):+
:X_i^{+ (i,2)}(z_1)S_i^{(i+1,1)}(z_2):\right).
\end{eqnarray}
Using the specialization (\ref{A2-12}), we conclude
$[X_i^+(z_1),S_i(z_2)]=0$ for $1\leq i \leq N-1$.
\\
$\bullet$~Proof of (\ref{thm:screening2}) 
for $1\leq i \neq j \leq N$.
\\
First we show (\ref{thm:screening2}) for $1\leq i \leq N-1, 
1\leq j\leq i-1$.
The commutators vanish,
$[X_i^{+ (l,s)}(z_1), S_j^{(m,t)}(z_2)]=0$,
for the following condition.
\begin{eqnarray}
[(l,s),(m,t)]
\neq \left\{\begin{array}{cc}
\begin{array}{c}
~[(j,1),(i,2)], [(j+1,2),(i+1,1)]\\
~[(j,1),(i+1,2)], [(j,2),(i+1,1)]\\
~[(j,2),(i,2)], [(j+1,2),(i+1,2)]
\end{array}
&~~(1\leq i \leq N-1, 1\leq j \leq i-1)
\end{array}\right..
\end{eqnarray}
Hence, using the relations (\ref{A1-21}), (\ref{A1-22}), 
(\ref{A1-25}), (\ref{A1-26}), (\ref{A1-27}), (\ref{A1-28}),
we have the following relation
for $1\leq i \leq N-1$ and $1\leq j \leq i-1$.
\begin{eqnarray}
&&[X_i^+(z_1),S_j(z_2)](q-q^{-1})z_1z_2\nonumber\\
&=&\delta\left(\frac{q^{N-i-j-1}z_2}{z_1}\right)
\left(
q:X_i^{+ (j,1)}(z_1)S_j^{(i,2)}(z_2):
-q:X_i^{+ (j+1,2)}(z_1)S_j^{(i+1,1)}(z_2):
\right.\nonumber\\
&&
-:X_i^{+ (j,1)}(z_1)S_j^{(i+1,2)}(z_2):
+:X_i^{+ (j,2)}(z_1)S_j^{(i+1,1)}(z_2):\nonumber\\
&&\left.
-q:X_i^{+ (j,2)}(z_1)S_j^{(i,2)}(z_2):
+q:X_i^{+ (j+1,2)}(z_1)S_j^{(i+1,2)}(z_2):\right).
\end{eqnarray}
Using the specializations (\ref{A2-13}), (\ref{A2-15}), (\ref{A2-16}), 
we conclude
$[X_i^+(z_1),S_j(z_2)]=0$ for $1\leq i \leq N-1, 1\leq j \leq i-1$.
Next we show 
(\ref{thm:screening2}) for $i=N$ and 
$1\leq j\leq N-1$.
The commutators vanish,
$[X_i^{+ (l,s)}(z_1), S_j^{(m,t)}(z_2)]=0$,
for the following condition.
\begin{eqnarray}
[(l,s),(m,t)]
\neq [(j,0),(N,2)], [(j+1,0),(N+1,0)]~~(1\leq j \leq N-1).
\end{eqnarray}
Hence, using the relations (\ref{A1-23}) and (\ref{A1-24}), we have
\begin{eqnarray}
&&[X_N^+(z_1),S_j(z_2)] q^{-N+1}z_2\nonumber\\
&=&\delta\left(\frac{q^{-j-1}z_2}{z_1}\right)
\left(:X_N^{+ (j,0)}(z_1)S_j^{(N,2)}(z_2):-
:X_N^{+ (j+1,0)}(z_1)S_j^{(N+1,0)}(z_2):\right).
\end{eqnarray}
Using the specialization (\ref{A2-14}),
we conclude $[X_N^+(z_1),S_j(z_2)]=0$ for $1\leq j \leq N-1$.
Now we have shown the commutation relation (\ref{thm:screening2}).\\
$\bullet$~Proof of (\ref{thm:screening1}). 
This is a direct consequence of the relation,
\begin{eqnarray}
[\Psi_i^\pm(z_1),S_j^{(l,s)}(z_2)]=0,
\end{eqnarray}
for every $i,j$, and $(l,s)$.\\
$\bullet$~
Proof of (\ref{thm:screening4}).
This is a direct consequence of the relation,
\begin{eqnarray}
\left[u_1-u_2+\frac{A_{i,j}}{2}\right]_{k+N-1}
S_i^{(l,s)}(z_1)S_j^{(m,t)}(z_2)=
\left[u_2-u_1+\frac{A_{i,j}}{2}\right]_{k+N-1}
S_j^{(m,t)}(z_2)S_i^{(l,s)}(z_1),
\end{eqnarray}
for every $(i,j)$ and $[(l,s),(m,t)]$.

\section{Vertex operator}

In this section we propose 
bosonizations of the vertex operators
for the quantum superalgebra $U_q(\widehat{sl}(N|1))$
\cite{Frenkel-Reshetikhin}.
We check that the vertex operators are the intertwiners among
the Fock-Wakimoto module and the typical representation
for small rank $N \leq 4$.

\subsection{Level-zero representation}

We discuss
level-zero representation of $U_q(\widehat{sl}(3|1))$
in this section 
(resp. $U_q(\widehat{sl}(4|1))$ in appendix \ref{appendixB}), 
that we will use for the investigation of
the vertex operator.
Let $V_\alpha$ be the one parameter family of
the $2^N$-dimensional typical representation
of $U_q(sl(N|1))$ \cite{Kac3, Palev-Tolstoy}.
In the case of $U_q(sl(3|1))$,
we choose the basis $\{v_j\}_{1\leq j \leq 8}$ of $V_\alpha$
and assign them the ${\bf Z}_2$-gradings as following.
\begin{eqnarray}
|v_1|=|v_5|=|v_6|=|v_7|=0,~~
|v_2|=|v_3|=|v_4|=|v_8|=1.
\end{eqnarray}
In the homogeneous gradation, 
the evaluation representation $V_{\alpha,z}$ of 
$U_q(\widehat{sl}(3|1))$ is given by
\begin{eqnarray}
h_1&=&E_{3,3}-E_{4,4}+E_{5,5}-E_{6,6},\\
h_2&=&E_{2,2}-E_{3,3}+E_{6,6}-E_{7,7},\\
h_3&=&\alpha(E_{1,1}+E_{2,2})+(\alpha+1)(E_{3,3}+E_{4,4}+E_{5,5}+E_{6,6})+
(\alpha+2)(E_{7,7}+E_{8,8}),\\
e_1&=&E_{3,4}+E_{5,6},\\
e_2&=&E_{2,3}+E_{6,7},\\
e_3&=&\sqrt{[\alpha]}E_{1,2}-\sqrt{[\alpha+1]}(E_{3,5}+E_{4,6})
+\sqrt{[\alpha+2]}E_{7,8},\\
f_1&=&E_{4,3}+E_{6,5},\\
f_2&=&E_{3,2}+E_{7,6},\\
f_3&=&\sqrt{[\alpha]}E_{2,1}-\sqrt{[\alpha+1]}(E_{5,3}+E_{6,4})+
\sqrt{[\alpha+2]}E_{8,7},\\
h_0&=&-\alpha(E_{1,1}+E_{4,4})
-(\alpha+1)(E_{2,2}+E_{3,3}+E_{6,6}+E_{7,7})
-(\alpha+2)(E_{5,5}+E_{8,8}),\\
e_0&=&
-z(\sqrt{[\alpha]}E_{4,1}
-\sqrt{[\alpha+1]}(E_{6,2}+E_{7,3})
+\sqrt{[\alpha+2]}E_{8,5}),\\
f_0&=&
z^{-1}(\sqrt{[\alpha]}E_{1,4}
-\sqrt{[\alpha+1]}(E_{2,6}+E_{3,7})
+\sqrt{[\alpha+2]}E_{5,8}).
\end{eqnarray}
We set the dual representation $V_{\alpha,z}^{* S}$ of
$U_q(\widehat{sl}(3|1))$ by
\begin{eqnarray}
\pi_{V_{\alpha,z}^{* S}}(a)=\left(\pi_{V_{\alpha,z}}
(S(a))\right)^{st}~~~
{\rm for}~~a \in U_q(\widehat{sl}(3|1)),
\end{eqnarray}
where we have used the antipode $S$ and 
have introduced the supertransposition "$st$" by
\begin{eqnarray}
(E_{i,j})^{st}=(-1)^{|v_i|(|v_i|+|v_j|)}E_{j,i}.
\end{eqnarray}
We have chosen the dual basis $\{v_j^*\}_{1\leq j \leq 8}$ 
of $V_\alpha^{* S}$
and assign them the ${\bf Z}_2$-gradings as following.
\begin{eqnarray}
|v_1^*|=|v_5^*|=|v_6^*|=|v_7^*|=0,~~
|v_2^*|=|v_3^*|=|v_4^*|=|v_8^*|=1.
\end{eqnarray}
In the homogeneous gradation, 
the evaluation representation $V_{\alpha,z}^{* S}$ of 
$U_q(\widehat{sl}(3|1))$ is given by
\begin{eqnarray}
h_1&=&-E_{3,3}+E_{4,4}-E_{5,5}+E_{6,6},\\
h_2&=&-E_{2,2}+E_{3,3}-E_{6,6}+E_{7,7},\\
h_3&=&-\alpha(E_{1,1}+E_{2,2})-(\alpha+1)(E_{3,3}+E_{4,4}+E_{5,5}+E_{6,6})
-(\alpha+2)(E_{7,7}+E_{8,8}),\\
e_1&=&-q^{-1}(E_{4,3}+E_{6,5}),\\
e_2&=&-q^{-1}(E_{3,2}+E_{7,6}),\\
e_3&=&-(
\sqrt{[\alpha]}q^{-\alpha}E_{2,1}
+\sqrt{[\alpha+1]}q^{-\alpha-1}(E_{5,3}+E_{6,4})
+\sqrt{[\alpha+2]}q^{-\alpha-2}E_{8,7}),\\
f_1&=&-q(E_{3,4}+E_{5,6}),\\
f_2&=&-q(E_{2,3}+E_{6,7}),\\
f_3&=&\sqrt{[\alpha]}q^{\alpha}E_{1,2}
+\sqrt{[\alpha+1]}q^{\alpha+1}(E_{3,5}+E_{4,6})+
\sqrt{[\alpha+2]}q^{\alpha+2}E_{7,8},\\
h_0&=&\alpha(E_{1,1}+E_{4,4})
+(\alpha+1)(E_{2,2}+E_{3,3}+E_{6,6}+E_{7,7})
+(\alpha+2)(E_{5,5}+E_{8,8}),
\\
e_0&=&
-z(\sqrt{[\alpha]}q^{\alpha}E_{1,4}
+\sqrt{[\alpha+1]}q^{\alpha+1}(E_{2,6}+E_{3,7})
+\sqrt{[\alpha+2]}q^{\alpha+2}E_{5,8}),\\
f_0&=&
-z^{-1}(\sqrt{[\alpha]}q^{-\alpha}E_{4,1}
+\sqrt{[\alpha+1]}q^{-\alpha-1}(E_{6,2}+E_{7,3})
+\sqrt{[\alpha+2]}q^{-\alpha-2}E_{8,5}).
\end{eqnarray}

We give the level-zero realization of the Drinfeld generators.

\begin{prop}~~~
On $V_{\alpha,z}$,
the Drinfeld generators of $U_q(\widehat{sl}(3|1))$
are realized by
\begin{eqnarray}
h_{1,m}&=&\frac{[m]}{m}(q^{\alpha+2}z)^m
(q^{-m}E_{3,3}-q^mE_{4,4}+q^{-m}E_{5,5}-q^mE_{6,6}),\\
h_{2,m}&=&
\frac{[m]}{m}(q^{\alpha+2}z)^m
(q^{-2m}E_{2,2}-E_{3,3}+E_{6,6}-q^{2m}E_{7,7}),\\
h_{3,m}&=&\frac{1}{m}z^m
([\alpha m](E_{1,1}+E_{2,2})+[(\alpha+1)m] q^m
(E_{3,3}+E_{4,4}+E_{5,5}+E_{6,6})\nonumber
\\
&&+[(\alpha+2)m] q^{2m} (E_{7,7}+E_{8,8})),\\
x_{1,n}^+&=&(q^{\alpha+2}z)^n (E_{3,4}+E_{5,6}),\\
x_{2,n}^+&=&(q^{\alpha+2}z)^n (q^{-n}E_{2,3}+q^nE_{6,7}),\\
x_{3,n}^+&=&(q^{\alpha+2} z)^n
(\sqrt{[\alpha]}q^{-2n}E_{1,2}-\sqrt{[\alpha+1]}
(E_{3,5}+E_{4,6})+\sqrt{[\alpha+2]}q^{2n}E_{7,8}),\\
x_{1,n}^-&=&
(q^{\alpha+2}z)^n (E_{4,3}+E_{6,5}),\\
x_{2,n}^-&=&(q^{\alpha+2}z)^n 
(q^{-n}E_{3,2}+q^nE_{7,6}),\\
x_{3,n}^-&=&(q^{\alpha+2} z)^n
(\sqrt{[\alpha]}q^{-2n}
E_{2,1}-\sqrt{[\alpha+1]}
(E_{5,3}+E_{6,4})+\sqrt{[\alpha+2]}q^{2n}E_{8,7}).
\end{eqnarray}
On $V_{\alpha,z}^{* S}$,
the Drinfeld generators of $U_q(\widehat{sl}(3|1))$
are realized by
\begin{eqnarray}
h_{1,m}&=&\frac{[m]_q}{m}(q^{-\alpha-2}z)^m
(-q^{m}E_{3,3}+q^{-m}E_{4,4}-q^{m}E_{5,5}+q^{-m}E_{6,6}),\\
h_{2,m}&=&
\frac{[m]_q}{m}(q^{-\alpha-2}z)^m
(-q^{2m}E_{2,2}+E_{3,3}-E_{6,6}+q^{-2m}E_{7,7}),\\
h_{3,m}&=&\frac{-1}{m}z^m
([\alpha m](E_{1,1}+E_{2,2})
+[(\alpha+1)m] q^{-m}
(E_{3,3}+E_{4,4}+E_{5,5}+E_{6,6})\nonumber
\\
&&+[(\alpha+2)m] q^{-2m} (E_{7,7}+E_{8,8})),\\
x_{1,n}^+&=&-q^{-1}(q^{-\alpha-2}z)^n (E_{4,3}+E_{6,5}),\\
x_{2,n}^+&=&-q^{-1}
(q^{-\alpha-2}z)^n (q^{n}E_{3,2}+q^{-n}E_{7,6}),\\
x_{3,n}^+&=&-
(q^{-\alpha-2} z)^n
(\sqrt{[\alpha]}q^{-\alpha+2n}E_{2,1}
+\sqrt{[\alpha+1]}q^{-\alpha-1}
(E_{5,3}+E_{6,4})\nonumber\\
&&+\sqrt{[\alpha+2]}q^{-\alpha-2-2n}E_{8,7}),
\\
x_{1,n}^-&=&-q
(q^{-\alpha}z)^n (E_{3,4}+E_{5,6}),\\
x_{2,n}^-&=&-q
(q^{-\alpha-2}z)^n 
(q^{n}E_{2,3}+q^{-n}E_{6,7}),\\
x_{3,n}^-&=&(q^{-\alpha-2} z)^n
(\sqrt{[\alpha]}
q^{\alpha+2n}E_{1,2}+\sqrt{[\alpha+1]}q^{\alpha+1}
(E_{3,5}+E_{4,6})\nonumber\\
&&+\sqrt{[\alpha+2]}q^{\alpha+2-2n}E_{7,8}).
\end{eqnarray}
\end{prop}
In appendix \ref{appendixB},
we summarize the case of $U_q(\widehat{sl}(4|1))$.
The case of $U_q(\widehat{sl}(2|1))$ is summarized in \cite{Zhang-Gould}.

\subsection{Vertex operator}

Let ${\cal F}$ and ${\cal F}'$ be level $k$ highest weight 
$U_q(\widehat{sl}(N|1))$-modules.
Let $V_{\alpha}$ and $V_{\alpha}^{* S}$ be
$2^N$-dimensional typical representation with a parameters $\alpha$
\cite{Palev-Tolstoy}.
The representations $V_{\alpha}$ and $V_\alpha^{* S}$
are irreducible if and only if 
$\alpha \neq 0, -1, -2, \cdots, -N+1$.
Let $V_{\alpha, z}$ and $V_{\alpha,z}^{* S}$ be
the evaluation module and its dual
of the typical representation.
Consider the following intertwiners 
of $U_q(\widehat{sl}(N|1))$-module
\cite{Frenkel-Reshetikhin}.
\begin{eqnarray}
\Phi(z): {\cal F} \longrightarrow {\cal F}' \otimes V_{\alpha,z},~~~
\Phi^*(z): {\cal F} \longrightarrow {\cal F}' \otimes V_{\alpha,z}^{* S}.
\end{eqnarray}
They are intertwiners in the sense that for any 
$x \in U_q(\widehat{sl}(N|1))$,
\begin{eqnarray}
\Phi(z)\cdot x=\Delta(x) \cdot \Phi(z),~~~
\Phi^*(z)\cdot x=\Delta(x) \cdot \Phi^*(z).
\end{eqnarray}
We expand the intertwining operators.
\begin{eqnarray}
\Phi(z)=\sum_{j=1}^{2^N}\Phi_j(z)\otimes v_j,~~~
\Phi^*(z)=\sum_{j=1}^{2^N}\Phi_j^*(z)\otimes v_j^*.
\end{eqnarray}
We set the ${\bf Z}_2$-grading of 
the intertwiner be $|\Phi(z)|=|\Phi^*(z)|=0$.
In what follow we focus our attention on
rank $N \leq 4$ case.

\begin{prop}
\label{prop:vertex1}
~~~For $\alpha \neq 0,-1,-2$, 
the operator $\Phi(z)$ for $U_q(\widehat{sl}(3|1))$
is determined by the component
$\Phi_8(z)$. More explicitly, we have
\begin{eqnarray}
\Phi_3(z)&=&[\Phi_4(z),f_1]_q,~~~
\Phi_5(z)=
[\Phi_6(z),f_1]_q,\\
\Phi_2(z)&=&[\Phi_3(z),f_2]_q,~~~
\Phi_6(z)=[\Phi_7(z),f_2]_q,
\\
\Phi_1(z)&=&
\frac{1}{\sqrt{[\alpha]}}
[\Phi_2(z),f_3]_{q^{-\alpha}},~~~
\Phi_3(z)=
\frac{-1}{\sqrt{[\alpha+1]}}
[\Phi_5(z),f_3]_{q^{-\alpha-1}},
\\
\Phi_4(z)&=&
\frac{-1}{\sqrt{[\alpha+1]}}
[\Phi_6(z),f_3]_{q^{-\alpha-1}},~~~
\Phi_7(z)=
\frac{1}{\sqrt{[\alpha+2]}}
[\Phi_8(z),f_3]_{q^{-\alpha-2}}.
\end{eqnarray}
For $\alpha \neq 0,-1,-2$, the operator
$\Phi^*(z)$ for $U_q(\widehat{sl}(3|1))$
is determined by the component
$\Phi_1^*(z)$. 
More explicitly, we have
\begin{eqnarray}
\Phi_4^*(z)&=&[f_1,\Phi_3^*(z)]_{q^{-1}},~~~
\Phi_6^*(z)=[f_1,\Phi_5^*(z)]_{q^{-1}},\\
\Phi_3^*(z)&=&[f_2,\Phi_2^*(z)]_{q^{-1}},~~~
\Phi_7^*(z)=[f_2,\Phi_6^*(z)]_{q^{-1}},
\\
\Phi_2^*(z)&=&
\frac{1}{\sqrt{[\alpha]}}
[f_3,\Phi_1^*(z)]_{q^{-\alpha}},~~~
\Phi_5^*(z)=
\frac{-1}{\sqrt{[\alpha+1]}}
[f_3,\Phi_3^*(z)]_{q^{-\alpha-1}},
\\
\Phi_6^*(z)&=&
\frac{-1}{\sqrt{[\alpha+1]}}
[f_3,\Phi_4^*(z)]_{q^{-\alpha-1}},~~~
\Phi_8^*(z)=
\frac{1}{\sqrt{[\alpha+2]}}
[f_3,\Phi_7^*(z)]_{q^{-\alpha-2}}.
\end{eqnarray}
\end{prop}

\begin{prop}
\label{prop:vertex2}
~~~For $\alpha \neq 0,-1,-2,-3$, the operator
$\Phi(z)$ for $U_q(\widehat{sl}(4|1))$
is determined by the component
$\Phi_{16}(z)$. More explicitly, we have
\begin{eqnarray}
\Phi_4(z)&=&[\Phi_6(z),f_1]_q,~~~
\Phi_7(z)=
[\Phi_8(z),f_1]_q,\\
\Phi_9(z)&=&[\Phi_{10}(z),f_1]_q,~~~
\Phi_{11}(z)=[\Phi_{13}(z),f_1]_q,
\\
\Phi_3(z)&=&[\Phi_4(z),f_2]_q,~~~
\Phi_5(z)=
[\Phi_7(z),f_2]_q,\\
\Phi_9(z)&=&[\Phi_{10}(z),f_2]_q,~~~
\Phi_{11}(z)=[\Phi_{13}(z),f_2]_q,
\\
\Phi_2(z)&=&[\Phi_3(z),f_3]_q,~~~
\Phi_7(z)=[\Phi_9(z),f_3]_q,\\
\Phi_8(z)&=&[\Phi_{10}(z),f_3]_q,~~~
\Phi_{14}(z)=[\Phi_{15}(z),f_2]_q,
\\
\Phi_1(z)&=&
\frac{-1}{\sqrt{[\alpha]}}
[\Phi_2(z),f_4]_{q^{-\alpha}},~~~
\Phi_3(z)=
\frac{-1}{\sqrt{[\alpha+1]}}
[\Phi_5(z),f_4]_{q^{-\alpha-1}},
\\
\Phi_4(z)&=&
\frac{-1}{\sqrt{[\alpha+1]}}
[\Phi_7(z),f_4]_{q^{-\alpha-1}},~~~
\Phi_6(z)=
\frac{-1}{\sqrt{[\alpha+1]}}
[\Phi_8(z),f_4]_{q^{-\alpha-1}},
\\
\Phi_9(z)&=&
\frac{-1}{\sqrt{[\alpha+2]}}
[\Phi_{11}(z),f_4]_{q^{-\alpha-2}},~~~
\Phi_{10}(z)=
\frac{-1}{\sqrt{[\alpha+2]}}
[\Phi_{13}(z),f_4]_{q^{-\alpha-2}},
\\
\Phi_{12}(z)&=&
\frac{-1}{\sqrt{[\alpha+2]}}
[\Phi_{14}(z),f_4]_{q^{-\alpha-2}},~~~
\Phi_{15}(z)=
\frac{-1}{\sqrt{[\alpha+3]}}
[\Phi_{16}(z),f_4]_{q^{-\alpha-3}}.
\end{eqnarray}
For $\alpha \neq 0,-1,-2,-3$, the operator
$\Phi^*(z)$ for $U_q(\widehat{sl}(4|1))$
is determined by the component
$\Phi_1^*(z)$. More explicitly, we have
\begin{eqnarray}
\Phi_6^*(z)&=&[f_1,\Phi_4^*(z)]_{q^{-1}},~~~
\Phi_8^*(z)=
[f_1,\Phi_7^*(z)]_{q^{-1}},\\
\Phi_{10}^*(z)&=&[f_1,\Phi_{9}^*(z)]_{q^{-1}},~~~
\Phi_{13}^*(z)=[f_1,\Phi_{11}^*(z)]_{q^{-1}},
\\
\Phi_4^*(z)&=&[f_2,\Phi_3^*(z)]_{q^{-1}},~~~
\Phi_7^*(z)=
[f_2,\Phi_5^*(z)]_{q^{-1}},\\
\Phi_{10}^*(z)&=&[f_2,\Phi_{9}^*(z)]_{q^{-1}},~~~
\Phi_{13}^*(z)=[f_2,\Phi_{11}^*(z)]_{q^{-1}},
\\
\Phi_3^*(z)&=&[f_3,\Phi_2^*(z)]_{q^{-1}},~~~
\Phi_9^*(z)=[f_3,\Phi_7^*(z)]_{q^{-1}},\\
\Phi_{10}^*(z)&=&[f_3,\Phi_{8}^*(z)]_{q^{-1}},~~~
\Phi_{15}^*(z)=[f_3,\Phi_{14}^*(z)]_{q^{-1}},
\\
\Phi_2^*(z)&=&
\frac{-1}{\sqrt{[\alpha]}}
[f_4,\Phi_1^*(z)]_{q^{-\alpha}},~~~
\Phi_5^*(z)=
\frac{1}{\sqrt{[\alpha+1]}}
[f_4,\Phi_3^*(z)]_{q^{-\alpha-1}},
\\
\Phi_7^*(z)&=&
\frac{1}{\sqrt{[\alpha+1]}}
[f_4,\Phi_4^*(z)]_{q^{-\alpha-1}},~~~
\Phi_8^*(z)=
\frac{1}{\sqrt{[\alpha+1]}}
[f_4,\Phi_6^*(z)]_{q^{-\alpha-1}},
\\
\Phi_{11}^*(z)&=&
\frac{-1}{\sqrt{[\alpha+2]}}
[f_4,\Phi_{9}^*(z)]_{q^{-\alpha-2}},~~~
\Phi_{13}^*(z)=
\frac{-1}{\sqrt{[\alpha+2]}}
[f_4,\Phi_{10}^*(z)]_{q^{-\alpha-2}},
\\
\Phi_{14}^*(z)&=&
\frac{-1}{\sqrt{[\alpha+2]}}
[f_4,\Phi_{12}^*(z)]_{q^{-\alpha-2}},~~~
\Phi_{16}^*(z)=
\frac{-1}{\sqrt{[\alpha+3]}}
[f_4,\Phi_{15}^*(z)]_{q^{-\alpha-3}}.
\end{eqnarray}
\end{prop}
The case of $U_q(\widehat{sl}(2|1))$ is summarized
in \cite{Zhang-Gould}.
Next we determine the relations between the components $\Phi_{2^N}(z)$,
$\Phi_1^*(z)$ and the Drinfeld generators.
We use the coproduct
\begin{eqnarray}
&&\Delta(h_i)=h_i \otimes 1+1\otimes h_i,\\
&&\Delta(h_{i,m})=h_{i,m}\otimes q^{\frac{cm}{2}}
+q^{\frac{3cm}{2}}\otimes h_{i,m}~~~(m>0),\\
&&\Delta(h_{i,-m})=h_{i,-m}\otimes q^{-\frac{3cm}{2}}
+q^{-\frac{cm}{2}}\otimes h_{i,-m}~~~(m>0).
\end{eqnarray}

\begin{prop}~~~
The component $\Phi_{2^N}(z)$ 
associated with $U_q(\widehat{sl}(N|1))$
satisfy
\begin{eqnarray}
&&[h_i, \Phi_{2^N}(z)]=-\delta_{i,N}(\alpha+N-1)\Phi_{2^N}(z)
~~~(1\leq i \leq N),\\
&&[h_{i,m},\Phi_{2^N}(z)]=-\delta_{i,N}q^{(N-1+\frac{3k}{2})m}
\frac{[(\alpha+N-1)m]}{m}z^m \Phi_{2^N}(z)~~~(m>0, 1\leq i \leq N),\\
&&
[h_{i,-m},\Phi_{2^N}(z)]=-\delta_{i,N}q^{(-N+1-\frac{k}{2})m}
\frac{[(\alpha+N-1)m]}{m}z^{-m} \Phi_{2^N}(z)~~~(m>0, 1\leq i \leq N),\\
&&
[X_i^+(z_1),\Phi_{2^N}(z_2)]=0~~~(1\leq i \leq N).
\end{eqnarray}
The component $\Phi_{1}^*(z)$ 
associated with $U_q(\widehat{sl}(N|1))$
satisfy
\begin{eqnarray}
&&[h_i, \Phi_{1}^*(z)]=\delta_{i,N} \alpha \Phi_{1}^*(z)
~~~(1\leq i \leq N),\\
&&[h_{i,m},\Phi_{1}^*(z)]=\delta_{i,N}q^{\frac{3k}{2}m}
\frac{[\alpha m]}{m}z^m \Phi_{1}^*(z)~~~(m>0, 1\leq i \leq N),\\
&&
[h_{i,-m},\Phi_{1}^*(z)]=\delta_{i,N}q^{-\frac{k}{2}m}
\frac{[\alpha m]}{m}z^{-m} \Phi_{1}^*(z)~~~(m>0, 1\leq i \leq N),\\
&&
[X_i^+(z_1),\Phi_{1}^*(z_2)]=0~~~(1\leq i \leq N).
\end{eqnarray}
We have checked this proposition for rank $N=2,3,4$.
\end{prop}

In order to construct bosonizations of $\Phi_{2^N}(z)$
and $\Phi_1^*(z)$, we introduce a bosonic operator $\phi^{l_a}(z|\beta)$.

\begin{dfn}~~~For
$l_a=(l_a^1,l_a^2,\cdots,l_a^N) \in {\bf C}^N$ and 
$\beta \in {\bf C}$,
we set the bosonic operator $\phi^{ l_a}(z|\beta)$ by
\begin{eqnarray}
\phi^{ l_a}(z|\beta)=:\exp\left(\sum_{i,j=1}^N
\left(\frac{l_a^i}{k+N-1}
\frac{{\rm Min}(i,j)}{N-1}\frac{N-1-{\rm Max}(i,j)}{1}
a^j\right)(z|\beta)\right):.
\end{eqnarray}
We call the operator $\phi^{l_a}(z|\beta)$ the "elementary vertex operator".
\end{dfn}

\begin{prop}~~~The highest vector 
$|\lambda\rangle=|l_a,0,0\rangle$
of $U_q(\widehat{sl}(N|1))$ is 
created from the Fock vacuum $|0\rangle$
and $\phi^{l_a}(z|\beta)$.
\begin{eqnarray}
|\lambda \rangle=\lim_{z \to 0}
\phi^{ l_a}(z|\beta)|0\rangle.
\end{eqnarray}
Here $|\lambda\rangle$
is the highest weight vector of the highest weight
whose classical part $\bar{\lambda}=\sum_{i=1}^N
l_a^i \bar{\Lambda}_i$.
\end{prop}
The elementary vertex operators $\phi^{ l_a}(z|\beta)$
give rise to the following map.
\begin{eqnarray}
\phi^{ l_a}(z|\beta) : 
F(p_a) \longrightarrow F(p_a+l_a).
\end{eqnarray}
Using the inversion relation,
\begin{eqnarray}
\sum_{r=1}^N \frac{[A_{i,r}m]}{[m]}\frac{[{\rm Min}(r,j)m]
[(N-1-{\rm Max}(r,j))m]}
{[(N-1)m][m]}=\delta_{i,j},
\end{eqnarray}
we have the following proposition.

\begin{prop}~~~
The elementary vertex operators $\phi^{ l_a}(z| \beta)$
satisfy the following relations.
\begin{eqnarray}
&&[h_{i,m}, \phi^{ l_a}(z|\beta)]=
\frac{1}{m}[l_a^i m]
q^{-(\beta+\frac{N-1}{2})|m|}z^m \phi^{ l_a}(z|\beta)~~~~~(1\leq i \leq N),\\
&&[X_{i}^+(z_1), \phi^{ l_a}(z_2| \beta)]=0~~~~~(1\leq i \leq N),\\
&&(z_1-q^{l_a^i}z_2)X_{i}^-(z_1)\phi^{ l_a}
\left(z_2\left|
-\frac{k+N-1}{2}\right.\right)\nonumber\\
&&=
(q^{l_a^i}z_1-z_2)
\phi^{ l_a}\left(z_2
\left|-\frac{k+N-1}{2}\right.\right)X_{i}^-(z_1)~~~~~
(1\leq i \leq N).
\end{eqnarray}
\end{prop}

~

\begin{prop}~~~~
For $k=\alpha\neq 0,-1,-2,\cdots,-N+1$,
bosonizations
of the components $\Phi_{2N}(z)$ and $\Phi_1^*(z)$
associated with $U_q(\widehat{sl}(N|1))$ are given by
\begin{eqnarray}
\Phi_{2^N}(z)=\phi^{\hat{l}}\left(q^{k+N-1} z\left|-\frac{k+N-1}{2}\right.
\right),~~~
\Phi_1^*(z)=\phi^{\hat{l}^*}\left(q^{k}z
\left|-\frac{k+N-1}{2}\right.\right),
\end{eqnarray}
where 
we have set
$\hat{l}=-(0,\cdots,0,\alpha+N-1)$
and $\hat{l}^*=(0,\cdots,0,\alpha)$.
The other components $\Phi_j(z)$ and $\Phi_j^*(z)$ 
$(1\leq j \leq 2^N)$
are represented 
by multiple contour integrals of Drinfeld currents
(cf. propositions \ref{prop:vertex1} and \ref{prop:vertex2}).
We have checked this proposition
for $N=2,3,4$.
\end{prop}

These bosonizations of the vertex operators are determined from
the commutation relations with the superalgebra $U_q(\widehat{sl}(N|1))$.
The construction is completely independent of which infinite dimensional
modules the vertex operators intertwine.
In what follow we shall clarify on which
space these vertex operators act.
We balance the "background charge" of the vertex operators
by using the screening currents.
For $x=(x_1,x_2,\cdots,x_N)\in {\bf N}^N$,
we set the screening operator
\begin{eqnarray}
{\cal Q}^{(x)}=:Q_1^{x_1}Q_2^{x_2}\cdots Q_N^{x_N}:.
\end{eqnarray}
The screening operator ${\cal Q}^{(x)}$ give rise to the map,
\begin{eqnarray}
{\cal Q}^{(x)} : F(p_a) \longrightarrow F(p_a+\hat{x}).
\end{eqnarray}
Here $\hat{x}=(\hat{x}_1,\hat{x}_2,\cdots,\hat{x}_N)$, where
$\hat{x}_i=\sum_{j=1}^N A_{i,j}x_j$.
The ${\cal Q}^{(x)}$ commute with
the projection operator
$\eta_0 \xi_0$.
Hence we have the map on the Fock-Wakimoto module.
\begin{eqnarray}
{\cal Q}^{(x)} : {\cal F}(p_a) \longrightarrow {\cal F}(p_a+\hat{x}).
\end{eqnarray}

\begin{dfn}
For $k=\alpha\neq 0,-1,-2,\cdots,-N$, we set the bosonic operators
\begin{eqnarray}
\widetilde{\Phi}^{(x)}(z)=\sum_{j=1}^{2^N} \widetilde{\Phi}_j^{(x)}(z)
\otimes
v_j,~~~
\widetilde{\Phi}^{(y)*}(z)=\sum_{j=1}^{2^N} \widetilde{\Phi}_j^{(y)*}(z)
\otimes v_j^*.
\end{eqnarray}
Here we have set
\begin{eqnarray}
\widetilde{\Phi}_j^{(x)}(z)
&=&\eta_0\xi_0 \cdot
{\cal Q}^{(x)} \cdot \Phi_j(z) \cdot \eta_0\xi_0,\\
\widetilde{\Phi}_j^{(y)*}(z)&=&
\eta_0\xi_0 \cdot
{\cal Q}^{(y)} \cdot \Phi_j^*(z) \cdot \eta_0\xi_0.
\end{eqnarray}
where $x=(x_1,x_2,\cdots,x_N) \in {\bf N}^N$ and
$y=(y_1,y_2,\cdots,y_N) \in {\bf N}^N$.
We call the operators $\widetilde{\Phi}^{(x)}(z), \widetilde{\Phi}^{(y)*}(z)$ 
the "projected vertex operators" for $U_q(\widehat{sl}(N|1))$.
\end{dfn}

\begin{prop}~~~
For $k=\alpha\neq 0,-1,-2,\cdots,-N+1$,
the projected vertex operators 
$\widetilde{\Phi}^{(x)}(z)$ and $\widetilde{\Phi}^{(y)*}(z)$ are 
the intertwiners among the Fock-Wakimoto module and
the typical representation.
\begin{eqnarray}
\widetilde{\Phi}^{(x)}(z) &:& {\cal F}(p_a) 
\longrightarrow {\cal F}(p_a+\hat{l}+\hat{x}) \otimes V_{\alpha,z},\\
\widetilde{\Phi}^{(y)*}(z) &:& {\cal F}(p_a) 
\longrightarrow {\cal F}(p_a+\hat{l^*}+\hat{y}) \otimes V_{\alpha,z}^{* S}.
\end{eqnarray}
Here we have set $\hat{l}=-(0,\cdots,0,\alpha+N-1)$ and 
$\hat{l^*}=(0,\cdots,0,\alpha)$.
Here we have set
$\hat{x}=(\hat{x}_1,\hat{x}_2,\cdots,\hat{x}_N)$ and
$\hat{y}=(\hat{y}_1,\hat{y}_2,\cdots,\hat{y}_N)$ 
where $\hat{x}_i=\sum_{j=1}^N A_{i,j}x_j$
and
$\hat{y}_i=\sum_{j=1}^N A_{i,j}y_j$.
We have checked this proposition for rank $N=2,3,4$.
\end{prop}

\subsection{Correlation function}

In this section we discuss an application of
the projected vertex operators 
$\widetilde{\Phi}^{(x)}(z)$ and $\widetilde{\Phi}^{*(y)}(z)$.
We study non-vanishing property of
the correlation function
which is defined to be the trace of the vertex operators
over the Fock-Wakimoto module of $U_q(\widehat{sl}(N|1))$,
that is
\begin{eqnarray}
{\rm Tr}_{{\cal F}(l_a)}\left(q^{L_0}
\widetilde{\Phi}_{j_1}^{(x_{(1)})}(z_1)
\widetilde{\Phi}_{j_2}^{(x_{(2)})}(z_2)\cdots 
\widetilde{\Phi}_{j_n}^{(x_{(n)})}(z_n)
\right).
\end{eqnarray}
Here we propose 
the $q$-Virasoro operator $L_0$ for $k=\alpha \neq -N+1$
as follows.
\begin{eqnarray}
L_0&=&\frac{1}{2}\sum_{i,j=1}^N \sum_{m \in {\bf Z}}
:a_{-m}^i\frac{m^2}{[m][(k+N-1)m]}\frac{[{\rm Min}(i,j)m]
[(N-1-{\rm Max}(i,j))m]}
{[(N-1)m][m]}a_m^j:\nonumber\\
&&+\sum_{i,j=1}^N \frac{{\rm Min}(i,j)(N-1-{\rm Max}(i,j))}
{(k+N-1)(N-1)}a_0^j\nonumber\\
&&-\frac{1}{2}\sum_{1\leq i<j \leq N}\sum_{m \in {\bf Z}}
:b_{-m}^{i,j}\frac{m^2}{[m]^2}b_m^{i,j}:
+\frac{1}{2}\sum_{1\leq i<j \leq N}
\sum_{m \in {\bf Z}}:c_{-m}^{i,j}\frac{m^2}{[m]^2}c_m^{i,j}:
\nonumber\\
&&+\frac{1}{2}\sum_{1\leq i \leq N}\sum_{m \in {\bf Z}}
:b_{-m}^{i,N+1}\frac{m^2}{[m]^2}b_m^{i,N+1}:
+\frac{1}{2}\sum_{1\leq i \leq N}b_{0}^{i,N+1}.
\end{eqnarray}
The $L_0$ eigenvalue of $|l_a,0,0\rangle$ is 
$\frac{1}{2(k+N-1)}(\bar{\lambda}|\bar{\lambda}+2\bar{\rho})$,
where $\bar{\rho}=\sum_{i=1}^N \bar{\Lambda}_i$ and $\bar{\lambda}=
\sum_{i=1}^N l_a^i \bar{\Lambda}_i$.

\begin{prop}~~~The correlation function of
the vertex operators,
\begin{eqnarray}
{\rm Tr}_{{\cal F}(l_a)}\left(q^{L_0}
\widetilde{\Phi}_{j_1}^{(x_{(1)})}(z_1)
\widetilde{\Phi}_{j_2}^{(x_{(2)})}(z_2)\cdots 
\widetilde{\Phi}_{j_n}^{(x_{(n)})}(z_n)
\right)\neq 0,
\end{eqnarray}
if and only if $k=\alpha\neq 0,-1,-2,\cdots,-N+1$ and
$x_{(s)}=(x_{(s),1},x_{(s),2}\cdots,x_{(s),N}) \in {\bf N}^N$ $(1\leq s \leq n)$ 
satisfy the condition,
\begin{eqnarray}
\sum_{s=1}^n x_{(s),i}=\frac{n \cdot i}{N-1} \alpha
+n \cdot i~~~~~(1\leq i \leq N).
\label{cor:eq1}
\end{eqnarray}
\end{prop}
We note that there doesn't exist non-rational solution $k=\alpha \notin {\bf Q}$
of the relation (\ref{cor:eq1}).
Next, we consider the correlation function involving also dual vertex
operators.

\begin{prop}~The correlation function of
the vertex operators and the dual vertex operators,
\begin{eqnarray}
{\rm Tr}_{{\cal F}(l_a)}\left(q^{L_0}
\widetilde{\Phi}_{i_1}^{*(y_{(1)})}(w_1)
\widetilde{\Phi}_{i_2}^{*(y_{(2)})}(w_2)\cdots 
\widetilde{\Phi}_{i_m}^{*(y_{(m)})}(w_m)
\widetilde{\Phi}_{j_1}^{(x_{(1)})}(z_1)
\widetilde{\Phi}_{j_2}^{(x_{(2)})}(z_2)\cdots 
\widetilde{\Phi}_{j_n}^{(x_{(n)})}(z_n)
\right)\neq 0,
\end{eqnarray}
if and only if $k=\alpha\neq 0,-1,-2,\cdots,-N+1$, 
$x_{(s)}=(x_{(s),1},x_{(s),2},\cdots,x_{(s),N}) \in {\bf N}^N$ $(1\leq s \leq n)$ 
and $y_{(t)}=(y_{(t),1},y_{(t),2},\cdots,y_{(t),N})\in {\bf N}^N$ 
$(1\leq t \leq m)$ satisfy the condition
\begin{eqnarray}
\sum_{s=1}^n x_{(s),i}+\sum_{t=1}^m y_{(t),i}=
\frac{(n-m) i}{N-1}\alpha+n \cdot i~~~~~(1\leq i \leq N).
\label{cor:eq2}
\end{eqnarray}
\end{prop}
We note that there exist non-rational solutions $k=\alpha \notin {\bf Q}$
of the relation (\ref{cor:eq2}).
Upon $k=\alpha \notin {\bf Q}$,
the relation (\ref{cor:eq2}) is equivalent to
\begin{eqnarray}
m=n~~~{\rm and}~~~\sum_{s=1}^n (x_{(s),i}+y_{(s),i})=n \cdot i~~~~~(1\leq i \leq N).
\end{eqnarray}
We conclude that the screening operators $Q_i$
are needed to ensure non-vanishing property of
correlation functions.
In other words, we have to balance the "background charge"
of the vertex operators to construct non-zero correlation functions.
We can write down integral representations
of the correlation functions by using
bosonizations of the vertex operators
\cite{Bouwknegt-Ceresole-McCarthy-Nieuwenhuizen}.
It is open and nontrivial problem to deform these integral representations 
to convenient formulae for physical applications.

~\\


\subsection*{Acknowledgements}~~~
This work is supported by the Grant-in-Aid for
Scientific Research {\bf C} (21540228)
from Japan Society for Promotion of Science. 
The author would like to thank
Professor Etsuro Date, Professor Hiroyuki Yamane, Professor Masato Okado,
Professor Kenji Iohara, 
and Professor Hitoshi Konno for their interests to this work.

~\\

\begin{appendix}

\section{Normal orderings}
\label{appendixA}

In this appendix we summarize formulae of normal orderings.
In order to get the following delta-function formulae,
the following normal orderings are useful.
\begin{eqnarray}
&&\exp\left(a_+^i(q^{\frac{k+N-1}{2}}z_1)\right)
:\exp\left(-\left(\frac{1}{k+N-1}a^j\right)
\left(z_2 \left|\frac{k+N-1}{2}\right.\right)\right):\nonumber\\
&&=
::q^{-A_{i,j}}
\frac{\displaystyle (1-q^{A_{i,j}-k-N+1}z_2/z_1)}
{\displaystyle (1-q^{-A_{i,j}-k-N+1}z_2/z_1)},\nonumber
\\
&&
:\exp\left(-\left(\frac{1}{k+N-1}a^j\right)
\left(z_2 \left|\frac{k+N-1}{2}\right.\right)\right):
\exp\left(a_+^i(q^{\frac{k+N-1}{2}}z_1)\right)
=::1,\nonumber\\
&&
:\exp\left(-\left(\frac{1}{k+N-1}a^i\right)\left(z_1\left|\frac{k+N-1}{2}
\right.\right)\right):
\exp\left(a_-^j(q^{-\frac{k+N-1}{2}}z_2)\right)\nonumber\\
&&=::
\frac{\displaystyle(1-q^{A_{i,j}-k-N+1}z_2/z_1)}
{\displaystyle(1-q^{-A_{i,j}-k-N+1}z_2/z_1)},\nonumber\\
&&
\exp\left(a_-^i(q^{-\frac{k+N-1}{2}}z_2)\right)
:\exp\left(-\left(\frac{1}{k+N-1}a^j\right)
\left(z_1 \left|\frac{k+N-1}{2}\right.\right)\right):=
::q^{A_{i,j}}.\nonumber
\end{eqnarray}
In order to get the following specialization relations,
the following formula is useful.
\begin{eqnarray}
b^{i,j}(qz)-b^{i,j}(q^{-1}z)=
b_+^{i,j}(z)-b_-^{i,j}(z).\nonumber
\end{eqnarray}

\subsection{Delta-function}

\begin{eqnarray}
&&[X_i^{- (N,0)}(z_1),S_i^{(N+1,0)}(z_2)]\nonumber\\
&&=
\frac{-1}{(q-q^{-1})q^{k+N}z_1z_2}
\left(
\delta\left(\frac{q^{-N-k+1}z_2}{z_1}\right)-
\delta\left(\frac{q^{-N-k-1}z_2}{z_1}\right)
\right)::~~(1\leq i \leq N-1),
\label{A1-1}
\end{eqnarray}
\begin{eqnarray}
[X_i^{- (i,1)}(z_1),S_i^{(i+1,1)}(z_2)]
=(q-q^{-1})\delta\left(\frac{q^{N+k-1}z_2}{z_1}\right)::~~
(1\leq i \leq N-1),
\label{A1-2}
\end{eqnarray}
\begin{eqnarray}
[X_i^{- (N-1,1)}(z_1),S_i^{(N,2)}(z_2)]=(q^{-1}-q)
\delta\left(\frac{q^{-N-k-1}z_2}{z_1}\right)::~~(1\leq i \leq N-1),
\label{A1-3}
\end{eqnarray}
\begin{eqnarray}
[X_i^{- (l+1,2)}(z_1),S_i^{(l+2,1)}(z_2)]=(q-q^{-1})
\delta\left(\frac{q^{N-3-k-2l}z_2}{z_1}\right)::
~~(1\leq i \leq l \leq N-2),
\label{A1-4}
\end{eqnarray}
\begin{eqnarray}
[X_i^{- (l+1,1)}(z_1),S_i^{(l+2,2)}(z_2)]=(q^{-1}-q)
\delta\left(\frac{q^{N-3-k-2l}z_2}{z_1}\right)::~~
(1\leq i \leq l \leq N-3),
\label{A1-5}
\end{eqnarray}
\begin{eqnarray}
[X_i^{- (i,2)}(z_1),S_i^{(i+1,2)}(z_2)]=(q^{-1}-q)
\delta\left(\frac{q^{N-3-k-2i}z_2}{z_1}\right)::
~~(1\leq i \leq N-2),
\label{A1-6}
\end{eqnarray}
\begin{eqnarray}
[X_N^{- (N,1)}(z_1),S_N^{(N+1,0)}(z_2)]=\frac{1}{q^{-N-k+1}z_1}
\delta\left(\frac{q^{N+k-1}z_2}{z_1}\right)::,
\label{A1-7}
\end{eqnarray}
\begin{eqnarray}
[X_N^{- (N,2)}(z_1),S_N^{(N+1,0)}(z_2)]=
\frac{1}{q^{N+k-1}z_1}
\delta\left(\frac{q^{-N-k+1}z_2}{z_1}\right)::,
\label{A1-8}
\end{eqnarray}
\begin{eqnarray}
[X_i^{- (l,1)}(z_1),S_{i+1}^{(l+1,2)}(z_2)]=(q-q^{-1})
\delta\left(\frac{q^{N-k-2l-2}z_2}{z_1}\right)::~~
(1\leq i \leq N-2, i+1 \leq l \leq N-1),\label{A1-9}
\end{eqnarray}
\begin{eqnarray}
[X_i^{- (l,2)}(z_1),S_{i+1}^{(l+1,1)}(z_2)]=(q^{-1}-q)
\delta\left(\frac{q^{N-k-2l-2}z_2}{z_1}\right)::~~
(1\leq i \leq N-2, i+1 \leq l \leq N-1),\label{A1-10}
\end{eqnarray}
\begin{eqnarray}
~[X_i^{- (j,1)}(z_1),S_j^{(i+1,1)}(z_2)]=(q-q^{-1})
\delta\left(\frac{q^{N+k-i+j-1}z_2}{z_1}\right)::~~
(1\leq j<i \leq N-1),\label{A1-11}
\end{eqnarray}
\begin{eqnarray}
~[X_i^{- (j+1,1)}(z_1),S_j^{(i,1)}(z_2)]=(q^{-1}-q)
\delta\left(\frac{q^{N+k-i+j-1}z_2}{z_1}\right)::~~
(1\leq j<i \leq N),\label{A1-12}
\end{eqnarray}
\begin{eqnarray}
~[X_i^{- (j,2)}(z_1),S_j^{(i+1,1)}(z_2)]=(q-q^{-1})
\delta\left(\frac{q^{N+k-i+j-1}z_2}{z_1}\right)::~~
(1\leq j<i \leq N-1),
\label{A1-13}
\end{eqnarray}
\begin{eqnarray}
~[X_i^{- (j+1,1)}(z_1),S_j^{(i,2)}(z_2)]=(q^{-1}-q)
\delta\left(\frac{q^{N+k-i+j-1}z_2}{z_1}\right)::~~(1\leq j<i \leq N),
\label{A1-14}
\end{eqnarray}
\begin{eqnarray}
[X_i^{- (j,1)}(z_1),S_j^{(i,2)}(z_2)]=(1-q^{-2})
\delta\left(\frac{q^{N+k-i+j-1}z_2}{z_1}\right)::~~(1\leq j<i \leq N),
\label{A1-15}
\end{eqnarray}
\begin{eqnarray}
[X_i^{- (j,2)}(z_1),S_j^{(i,1)}(z_2)]=(q^{-2}-1)
\delta\left(\frac{q^{N+k-i+j-1}z_2}{z_1}\right)::~~(1\leq j<i \leq N),
\label{A1-16}
\end{eqnarray}
\begin{eqnarray}
[X_N^{- (j,1)}(z_1),S_j^{(N+1,0)}(z_2)]=\frac{1}{q^{-k-j+1}z_1}
\delta\left(\frac{q^{k+j-1}z_2}{z_1}\right)::~~(1\leq j \leq N-1),
\label{A1-17}
\end{eqnarray}
\begin{eqnarray}
[X_N^{- (j,2)}(z_1),S_j^{(N+1,0)}(z_2)]=\frac{1}{q^{-k-j+1}z_1}
\delta\left(\frac{q^{k+j-1}z_2}{z_1}\right)::~~(1\leq j \leq N-1),
\label{A1-18}
\end{eqnarray}
\begin{eqnarray}
[X_i^{+ (i,1)}(z_1), S_i^{(i+1,2)}(z_2)]=
(q-q^{-1})\delta
\left(\frac{q^{N-2i-1}z_2}{z_1}\right)
::~~(1\leq i \leq N-1),
\label{A1-19}
\end{eqnarray}
\begin{eqnarray}
[X_i^{+ (i,2)}(z_1), S_i^{(i+1,1)}(z_2)]=
(q^{-1}-q)\delta\left(\frac{q^{N-2i-1}z_2}{z_1}\right)
::~~(1\leq i \leq N-1),
\label{A1-20}
\end{eqnarray}
\begin{eqnarray}
[X_i^{+ (j,1)}(z_1),S_j^{(i,2)}(z_2)]=(1-q^2)
\delta\left(\frac{q^{N-i-j-1}z_2}{z_1}\right)
::~~
(1\leq i \leq N, 1\leq j \leq i-1),
\label{A1-21}
\end{eqnarray}
\begin{eqnarray}
[X_i^{+ (j+1,2)}(z_1),S_j^{(i+1,1)}(z_2)]=
(q^2-1)\delta\left(\frac{q^{N-i-j-1}z_2}{z_1}\right)
::~~
(1\leq i \leq N, 1\leq j \leq i-1),
\label{A1-22}
\end{eqnarray}
\begin{eqnarray}
[X_N^{+ (j,0)}(z_1), S_j^{(N,2)}(z_2)]&=&(1-q^2)
\delta\left(\frac{q^{-j-1}z_2}{z_1}\right)
::~~(1\leq j \leq N-1),
\label{A1-23}
\end{eqnarray}
\begin{eqnarray}
[X_N^{+ (j+1,0)}(z_1),S_j^{(N+1,0)}(z_2)]=
\frac{-1}{q^{j+1}z_1}\delta\left(\frac{q^{-j-1}z_2}{z_1}\right)
::~~(1\leq j \leq N-1),
\label{A1-24}
\end{eqnarray}
\begin{eqnarray}
[X_i^{+ (j,1)}(z_1), S_j^{(i+1,2)}(z_2)]&=&(q-q^{-1})\delta
\left(\frac{q^{N-i-j-1}z_2}{z_1}\right)::~~
(1\leq i \leq N-1, 1\leq j \leq i-1),
\label{A1-25}
\end{eqnarray}
\begin{eqnarray}
[X_i^{+ (j,2)}(z_1), S_j^{(i+1,1)}(z_2)]=
(q^{-1}-q)\delta
\left(\frac{q^{N-i-j-1}z_2}{z_1}\right)::~~
(1\leq i \leq N-1, 1\leq j \leq i-1),
\label{A1-26}
\end{eqnarray}
\begin{eqnarray}
[X_i^{+ (j,2)}(z_1), S_j^{(i,2)}(z_2)]=(1-q^2)
\delta\left(\frac{q^{N-i-j-1}z_2}{z_1}\right)::~~
(1\leq i \leq N-1, 1\leq j\leq i-1),
\label{A1-27}
\end{eqnarray}
\begin{eqnarray}
[X_i^{+ (j+1,2)}(z_1),S_j^{(i+1,2)}(z_2)]=(q^2-1)
\delta\left(\frac{q^{N-i-j-1}z_2}{z_1}\right)::
~~(1\leq i \leq N-1, 1\leq j \leq i-1),
\label{A1-28}
\end{eqnarray}

\subsection{Specialization}

\begin{eqnarray}
:X_i^{- (N,0)}(z)S_i^{(N+1,0)}(q^{N+k-1}z):=
:X_i^{- (i,1)}(z)S_i^{(i+1,1)}(q^{-N-k+1}z):\nonumber\\
=
:\exp\left(-\left(\frac{1}{k+N-1}a^i\right)\left(z\left|-\frac{k+N-1}{2}
\right.\right)\right):
~~~(1\leq i \leq N-1),
\label{A2-1}
\end{eqnarray}
\begin{eqnarray}
:X_N^{- (N,1)}(z)S_N^{(N+1,0)}(q^{-N-k+1}z):=
:X_N^{- (N,2)}(z)S_N^{(N+1,0)}(q^{N+k-1}z):\nonumber\\
=
:\exp\left(-\left(\frac{1}{k+N-1}a^i\right)\left(z\left|
-\frac{k+N-1}{2}\right.\right)\right):,
\label{A2-2}
\end{eqnarray}
\begin{eqnarray}
:X_i^{- (N,0)}(z)S_i^{(N+1,0)}(q^{N+k+1}z):
=:X_i^{- (N-1,1)}(z)S_i^{(N,2)}(q^{N+k+1}z):\nonumber\\
(1\leq i \leq N-1),
\label{A2-3}
\end{eqnarray}
\begin{eqnarray}
:X_i^{- (i+1,2)}(z)S_i^{(i+2,1)}(q^{-N+3+k+2i}z):
=:X_i^{- (i,2)}(z)S_i^{(i+1,2)}(q^{-N+3+k+2i}z):\nonumber
\\
(1\leq i \leq N-1),
\label{A2-4}
\end{eqnarray}
\begin{eqnarray}
:X_i^{- (l+1,1)}(z)S_i^{(l+2,2)}(q^{-N+3+k+2l}z):
=:X_i^{- (l+1,2)}(z)S_i^{(l+2,1)}(q^{-N+3+k+2l}z):\nonumber\\
(1\leq i \leq N-1, i+1 \leq l \leq N-2).\label{A2-5}
\end{eqnarray}
\begin{eqnarray}
:X_i^{- (l,1)}(z)S_{i+1}^{(l+1,2)}(q^{-N+k+2l+2}z):=
:X_i^{- (l,2)}(z)S_{i+1}^{(l+1,1)}(q^{-N+k+2l+2}z):,\nonumber\\
(1\leq i \leq N-2, i+1 \leq l \leq N-1), \label{A2-6}
\end{eqnarray}
\begin{eqnarray}
:X_i^{- (j,1)}(z)S_j^{(i+1,1)}(q^{-N-k+i-j+1}z):
=
:X_i^{- (j+1,1)}(z)S_j^{(i,1)}(q^{-N-k+i-j+1}z):\nonumber\\
(1\leq j<i \leq N-1),
\label{A2-7}
\end{eqnarray}
\begin{eqnarray}
:X_i^{- (j,2)}(z)S_j^{(i+1,1)}(q^{-N-k+i-j+1}z):
=
:X_i^{- (j+1,1)}(z)S_j^{(i,2)}(q^{-N-k+i-j+1}z):\nonumber\\
(1\leq j<i \leq N-1),
\label{A2-8}
\end{eqnarray}
\begin{eqnarray}
:X_i^{- (j,2)}(z)S_j^{(i,1)}(q^{-N-k+i-j+1}z):
=
:X_i^{- (j,1)}(z)S_j^{(i,2)}(q^{-N-k+i-j+1}z):\nonumber\\
(1\leq j<i \leq N),
\label{A2-9}
\end{eqnarray}
\begin{eqnarray}
:X_N^{- (j,1)}(z)S_j^{(N+1,0)}(q^{-k-j+1}z):
=
:X_N^{- (j+1,1)}(z)S_j^{(N,1)}(q^{-k-j+1}z):\nonumber\\
(1\leq j<i \leq N-1),
\label{A2-10}
\end{eqnarray}
\begin{eqnarray}
:X_N^{- (j,2)}(z)S_j^{(N+1,0)}(q^{-k-j+1}z):
=
:X_N^{- (j+1,1)}(z)S_j^{(N,2)}(q^{-k-j+1}z):\nonumber\\
(1\leq j<i \leq N-1),
\label{A2-11}
\end{eqnarray}
\begin{eqnarray}
:X_i^{+ (i,1)}(z)S_i^{(i+1,2)}(q^{-N+2i+1}z):=
:X_i^{+ (i,2)}(z)S_i^{(i+1,1)}(q^{-N+2i+1}z):\nonumber\\
(1\leq i \leq N-1),
\label{A2-12}
\end{eqnarray}
\begin{eqnarray}
:X_i^{+ (j,1)}(z)S_j^{(i,2)}(q^{-N+i+j+1}z):=
:X_i^{+ (j+1,2)}(z)S_j^{(i+1,1)}(q^{-N+i+j+1}z):\nonumber\\
(1\leq i \leq N-1, 1\leq j \leq i-1),
\label{A2-13}
\end{eqnarray}
\begin{eqnarray}
:X_N^{+ (j,0)}(z)S_j^{(N,2)}(q^{j+1}z):=
:X_N^{+ (j+1,0)}(z)S_j^{(N+1,0)}(q^{j+1}z):\nonumber\\
(1\leq j \leq N-1), 
\label{A2-14}
\end{eqnarray}
\begin{eqnarray}
:X_i^{+ (j,1)}(z)S_j^{(i+1,2)}(q^{-N+i+j+1}z):=
:X_i^{+ (j,2)}(z)S_j^{(i+1,1)}(q^{-N+i+j+1}z):
\nonumber\\
(1\leq i \leq N-1, 1\leq j \leq i-1),
\label{A2-15}
\end{eqnarray}
\begin{eqnarray}
:X_i^{+ (j,2)}(z)S_j^{(i,2)}(q^{-N+i+j+1}z):=
:X_i^{+ (j+1,2)}(z)S_j^{(i+1,2)}(q^{-N+i+j+1}z):\nonumber\\
(1\leq i \leq N-1, 1\leq j \leq i-1),
\label{A2-16}
\end{eqnarray}

\section{Level-zero representation of $U_q(\widehat{sl}(4|1))$}
\label{appendixB}

In this appendix we summarize
the level-zero representation of $U_q(\widehat{sl}(4|1))$.
Let $V_\alpha$ be the one parameter family of the $16(=2^4)$-dimensional
typical representation of $U_q(sl(4|1))$
\cite{Kac3, Palev-Tolstoy}.
In the case of $U_q(sl(4|1))$,
we choose the basis $\{v_j\}_{1\leq j \leq 16}$ 
of $V_\alpha$
and assign them the ${\bf Z}_2$-gradings as following.
\begin{eqnarray}
&&|v_1|=|v_5|=|v_7|=|v_8|=|v_9|=|v_{10}|=|v_{12}|=|v_{16}|=0,\nonumber
\\
&&|v_2|=|v_3|=|v_4|=|v_6|=|v_{11}|=|v_{13}|=|v_{14}|=|v_{15}|=1.
\end{eqnarray}
In the homogeneous gradation, 
the evaluation representation $V_{\alpha,z}$ of 
$U_q(\widehat{sl}(4|1))$ is given by
\begin{eqnarray}
h_1&=&E_{4,4}-E_{6,6}+E_{7,7}-E_{8,8}+E_{9,9}-E_{10,10}+E_{11,11}-E_{13,13},
\\
h_2&=&E_{3,3}-E_{4,4}+E_{5,5}-E_{7,7}+E_{10,10}-E_{12,12}+E_{13,13}-E_{14,14},
\\
h_3&=&E_{2,2}-E_{3,3}+E_{7,7}-E_{9,9}+E_{8,8}-E_{10,10}+E_{14,14}-E_{15,15},
\\
h_4&=& \alpha \sum_{j=1}^2 E_{j,j}
+(\alpha+1)\sum_{j=3}^8 E_{j,j}+(\alpha+2)
\sum_{j=9}^{14}E_{j,j}
+(\alpha+3)\sum_{j=15}^{16}E_{j,j},\\
e_1&=&E_{4,6}+E_{7,8}+E_{9,10}+E_{11,13},\\
e_2&=&E_{3,4}+E_{5,7}+E_{10,12}+E_{13,14},\\
e_3&=&E_{2,3}+E_{7,9}+E_{8,10}+E_{14,15},\\
e_4&=&-\sqrt{[\alpha]}E_{1,2}+\sqrt{[\alpha+1]}(E_{3,5}+E_{4,7}+E_{6,8})\nonumber\\
&&-\sqrt{[\alpha+2]}(E_{9,11}+E_{10,13}+E_{12,14})
+\sqrt{[\alpha+3]}E_{15,16},\\
f_1&=&E_{6,4}+E_{8,7}+E_{10,9}+E_{13,11},\\
f_2&=&E_{4,3}+E_{7,5}+E_{12,10}+E_{14,13},\\
f_3&=&E_{3,2}+E_{9,7}+E_{10,8}+E_{15,14},\\
f_4&=&-\sqrt{[\alpha]}E_{2,1}+\sqrt{[\alpha+1]}
(E_{5,3}+E_{7,4}+E_{8,6})\nonumber\\
&&-\sqrt{[\alpha+2]}(E_{11,9}+E_{13,10}+E_{14,12})+
\sqrt{[\alpha+3]}E_{16,15},\\
h_0&=&-\alpha(E_{1,1}+E_{6,6})
-(\alpha+1)(E_{2,2}+E_{3,3}+E_{4,4}+E_{8,8}
+E_{10,10}+E_{12,12})\nonumber\\
&&-(\alpha+2)(E_{5,5}+E_{7,7}+E_{9,9}
+E_{13,13}+E_{14,14}+E_{15,15})-(\alpha+3)(E_{11,11}+E_{16,16}),\\
e_0&=&z(\sqrt{[\alpha]}E_{6,1}-\sqrt{[\alpha+1]}
(E_{8,2}+E_{10,3}+E_{12,4})\nonumber\\
&&+\sqrt{[\alpha+2]}(E_{13,5}+E_{14,7}+E_{15,9})-
\sqrt{[\alpha+3]}E_{16,11}),\\
f_0&=&-z^{-1}(\sqrt{[\alpha]}E_{1,6}-\sqrt{[\alpha+1]}
(E_{2,8}+E_{3,10}+E_{4,12})\nonumber\\
&&+\sqrt{[\alpha+2]}(E_{5,13}+E_{7,14}+E_{9,15})-
\sqrt{[\alpha+3]}E_{11,16}).
\end{eqnarray}
We choose the dual basis $\{v_j^*\}_{1\leq j \leq 16}$ 
of $V_\alpha^{*,S}$
and assign them the ${\bf Z}_2$-gradings as following.
\begin{eqnarray}
&&|v_1^*|=|v_5^*|=|v_7^*|=|v_8^*|=|v_9^*|=|v_{10}^*|=
|v_{12}^*|=|v_{16}^*|=0,\nonumber
\\
&&|v_2^*|=|v_3^*|=|v_4^*|=|v_6^*|=|v_{11}^*|=
|v_{13}^*|=|v_{14}^*|=|v_{15}^*|=1.
\end{eqnarray}
In the homogeneous gradation, 
the evaluation representation $V_{\alpha,z}^{* S}$ of 
$U_q(\widehat{sl}(4|1))$ is given by
\begin{eqnarray}
h_1&=&
-E_{4,4}+E_{6,6}-E_{7,7}+E_{8,8}-E_{9,9}
+E_{10,10}-E_{11,11}+E_{13,13},
\\
h_2&=&-E_{3,3}+E_{4,4}-E_{5,5}+E_{7,7}
-E_{10,10}+E_{12,12}-E_{13,13}+E_{14,14},
\\
h_3&=&-E_{2,2}+E_{3,3}-E_{7,7}+E_{9,9}
-E_{8,8}+E_{10,10}-E_{14,14}+E_{15,15},
\\
h_4&=& -\alpha \sum_{j=1}^2E_{j,j}-(\alpha+1)
\sum_{j=3}^8 E_{j,j}-(\alpha+2)
\sum_{j=9}^{14}E_{j,j}
-(\alpha+3)\sum_{j=15}^{16}E_{j,j},
\\
e_1&=&-q^{-1}(E_{6,4}+E_{8,7}+E_{10,9}+E_{13,11}),\\
e_2&=&-q^{-1}(E_{4,3}+E_{7,5}+E_{12,10}+E_{14,13}),\\
e_3&=&-q^{-1}(E_{3,2}+E_{9,7}+E_{10,8}+E_{15,14}),\\
e_4&=&\sqrt{[\alpha]}q^{-\alpha}E_{2,1}
+\sqrt{[\alpha+1]}q^{-\alpha-1}(E_{5,3}+E_{7,4}+E_{8,6})
\nonumber\\
&&+\sqrt{[\alpha+2]}q^{-\alpha-2}
(E_{11,9}+E_{13,10}+E_{14,12})
+\sqrt{[\alpha+3]}q^{-\alpha-3}E_{16,15},
\\
f_1&=&-q(E_{4,6}+E_{7,8}+E_{9,10}+E_{11,13}),\\
f_2&=&-q(E_{3,4}+E_{5,7}+E_{10,12}+E_{13,14}),\\
f_3&=&-q(E_{2,3}+E_{7,9}+E_{8,10}+E_{14,15}),\\
f_4&=&-
\sqrt{[\alpha]}q^{\alpha}E_{1,2}
-\sqrt{[\alpha+1]}q^{\alpha+1}(E_{3,5}+E_{4,7}+E_{6,8})
\nonumber\\
&&-\sqrt{[\alpha+2]}q^{\alpha+2}
(E_{9,11}+E_{10,13}+E_{12,14})
-\sqrt{[\alpha+3]}q^{\alpha+3}E_{15,16},
\\
h_0&=&\alpha(E_{1,1}+E_{6,6})
+(\alpha+1)(E_{2,2}+E_{3,3}+E_{4,4}+E_{8,8}
+E_{10,10}+E_{12,12})\nonumber\\
&&+(\alpha+2)(E_{5,5}+E_{7,7}+E_{9,9}
+E_{13,13}+E_{14,14}+E_{15,15})
+(\alpha+3)(E_{11,11}+E_{16,16}),\\
e_0&=&z(\sqrt{[\alpha]}q^{\alpha}E_{1,6}
+\sqrt{[\alpha+1]}q^{\alpha+1}
(E_{2,8}+E_{3,10}+E_{4,12})\nonumber\\
&&+\sqrt{[\alpha+2]}q^{\alpha+2}
(E_{5,13}+E_{7,14}+E_{9,15})+
\sqrt{[\alpha+3]}q^{\alpha+3}E_{11,16}),
\\
f_0&=&z^{-1}(\sqrt{[\alpha]}q^{-\alpha}E_{6,1}
+\sqrt{[\alpha+1]}q^{-\alpha-1}
(E_{8,2}+E_{10,3}+E_{12,4})\nonumber\\
&&+\sqrt{[\alpha+2]}q^{-\alpha-2}
(E_{13,5}+E_{14,7}+E_{15,9})+
\sqrt{[\alpha+3]}q^{-\alpha-3}E_{16,11}).
\end{eqnarray}

We give the level-zero realization of the Drinfeld generators.

\begin{prop}
~~~On $V_{\alpha,z}$, 
the Drinfeld generators of $U_q(\widehat{sl}(4|1))$
are given by
\begin{eqnarray}
h_{1,m}&=&
\frac{[m]}{m}
(q^{\alpha+3}z)^m
(q^{-m}E_{4,4}-q^mE_{6,6}+
q^{-m}E_{7,7}-q^mE_{8,8}\nonumber\\
&+&q^{-m}E_{9,9}-q^mE_{10,10}+
q^{-m}E_{11,11}-q^mE_{13,13}),\\
h_{2,m}&=&
\frac{[m]}{m}(q^{\alpha+3}z)^m
(q^{-2m}E_{3,3}-E_{4,4}+q^{-2m}E_{5,5}-E_{7,7}\nonumber\\
&+&E_{10,10}-q^{2m}E_{12,12}+E_{13,13}-q^{2m}E_{14,14}),\\
h_{3,m}&=&
\frac{[m]}{m}(q^{\alpha+3}z)^m
(q^{-3m}E_{2,2}-q^{-m}E_{3,3}+q^{-m}E_{7,7}-q^{m}E_{9,9}\nonumber\\
&+&q^{-m}E_{8,8}-q^{m}E_{10,10}+q^{m}E_{14,14}-q^{3m}E_{15,15}),\\
h_{4,m}&=&\frac{1}{m}z^m
\left([\alpha m] \sum_{j=1}^2E_{j,j}
+[(\alpha+1)m] q^m \sum_{j=3}^8 E_{j,j}\right.\nonumber\\
&+&\left.[(\alpha+2)m] q^{2m} \sum_{j=9}^{14} E_{j,j}
+[(\alpha+3)m] q^{3m} \sum_{j=15}^{16}E_{j,j}\right),\\
x_{1,n}^+&=&
(q^{\alpha+3}z)^m(E_{4,6}+E_{7,8}+E_{9,10}+E_{11,13}),\\
x_{2,n}^+&=&
(q^{\alpha+3}z)^n(
q^{-3n}E_{3,4}+q^{-n}E_{5,7}+q^{n}E_{10,12}+q^{n}E_{13,14}),\\
x_{3,n}^+&=&
(q^{\alpha+3}z)^n(
q^{-2n}E_{2,3}+E_{7,9}+E_{8,10}+q^{2n}E_{14,15}),\\
x_{4,n}^+&=&(q^{\alpha+3}z)^n
(-\sqrt{[\alpha]}q^{-3n}E_{1,2}
+\sqrt{[\alpha+1]}q^{-n}(E_{3,5}+E_{4,7}+E_{6,8})
\nonumber\\
&-&\sqrt{[\alpha+2]_q}q^{n}(E_{9,11}+E_{10,13}+E_{12,14})
+\sqrt{[\alpha+3]}q^{3n}E_{15,16}),\\
x_{1,n}^-&=&
(q^{\alpha+3}z)^m(E_{6,4}+E_{8,7}+E_{10,9}+E_{13,11}),\\
x_{2,n}^-&=&
(q^{\alpha+3}z)^n(
q^{-n}E_{4,3}+q^{-n}E_{7,5}+q^{n}E_{12,10}+q^{n}E_{14,13}),\\
x_{3,n}^-&=&
(q^{\alpha+3}z)^n(
q^{-2n}E_{3,2}+E_{9,7}+E_{10,8}+q^{2n}E_{15,14}),\\
x_{4,n}^-&=&(q^{\alpha+3}z)^n
(-\sqrt{[\alpha]}q^{-3n}E_{2,1}
+\sqrt{[\alpha+1]}q^{-n}(E_{5,3}+E_{7,4}+E_{8,6})
\nonumber\\
&-&\sqrt{[\alpha+2]}q^{n}(E_{11,9}+E_{13,10}+E_{14,12})
+\sqrt{[\alpha+3]}q^{3n}E_{16,15}).
\end{eqnarray}
On $V_{\alpha,z}^{* S}$, the Drinfeld generators of 
$U_q(\widehat{sl}(4|1))$ are given by
\begin{eqnarray}
h_{1,m}&=&\frac{[m]}{m}(q^{-\alpha-3}z)^m
(-q^m E_{4,4}+q^{-m}E_{6,6}-q^mE_{7,7}+q^{-m}E_{8,8}
\nonumber\\
&&-q^m E_{9,9}+q^{-m}E_{10,10}-q^mE_{11,11}+q^{-m}E_{13,13}),\\
h_{2,m}
&=&\frac{[m]}{m}(q^{-\alpha-3}z)^m
(-q^{2m}E_{3,3}+E_{4,4}-q^{2m}E_{5,5}+E_{7,7}\nonumber\\
&&-E_{10,10}+q^{-2m}E_{12,12}-E_{13,13}+q^{-2m}E_{14,14}),\\
h_{3,m}
&=&\frac{[m]}{m}(q^{-\alpha-3}z)^m
(-q^{3m}E_{2,2}+q^m E_{3,3}-q^{m}E_{7,7}+q^{-m}E_{9,9}\nonumber\\
&&-q^mE_{8,8}+q^{-m}E_{10,10}-q^{-m}E_{14,14}
+q^{-3m}E_{15,15}),\\
h_{4,m}&=&\frac{-1}{m}z^m
([\alpha m]\sum_{j=1}^2 E_{j,j}+
[(\alpha+1)m] q^{-m} \sum_{j=3}^8E_{j,j}\nonumber\\
&&+[(\alpha+2)m]q^{-2m}\sum_{j=9}^{14}E_{j,j}
+[(\alpha+3)m]q^{-3m}\sum_{j=15}^{16}E_{j,j}),
\\
x_{1,n}^+&=&-q^{-1}(q^{-\alpha-3}z)^n
(E_{6,4}+E_{8,7}+E_{10,9}+E_{13,11}),\\
x_{2,n}^+&=&-q^{-1}
(q^{-\alpha-3}z)^n
(q^n E_{4,3}+q^n E_{7,5}+q^{-n}E_{12,10}
+q^{-n}E_{14,13}),\\
x_{3,n}^+&=&-q^{-1}
(q^{-\alpha-3}z)^n
(q^{2n} E_{3,2}+E_{9,7}+E_{10,8}+q^{-2n}E_{15,14}),\\
x_{4,n}^+&=&
(q^{-\alpha-3}z)^n
(\sqrt{[\alpha]}q^{-\alpha+3n}E_{2,1}
+\sqrt{[\alpha+1]}q^{-\alpha-1+n}(E_{5,3}+E_{7,4}+E_{8,6})
\nonumber\\
&&+\sqrt{[\alpha+2]}q^{-\alpha-2-n}
(E_{11,9}+E_{13,10}+E_{14,12})+
\sqrt{[\alpha+3]}q^{-\alpha-3-3n}E_{16,15}),
\\
x_{1,n}^-&=&-q(q^{-\alpha-3}z)^n
(E_{4,6}+E_{7,8}+E_{9,10}+E_{11,13}),\\
x_{2,n}^-&=&-q
(q^{-\alpha-3}z)^n
(q^n E_{3,4}+q^n E_{5,7}+q^{-n}E_{10,12}
+q^{-n}E_{13,14}),\\
x_{3,n}^-&=&-q
(q^{-\alpha-3}z)^n
(q^{2n} E_{2,3}+E_{7,9}+E_{8,10}+q^{-2n}E_{14,15}),\\
x_{4,n}^-&=&
-(q^{-\alpha-3}z)^n
(\sqrt{[\alpha]}q^{\alpha+3n}E_{1,2}
+\sqrt{[\alpha+1]}q^{\alpha+1+n}(E_{3,5}+E_{4,7}+E_{6,8})
\nonumber\\
&&+\sqrt{[\alpha+2]}q^{\alpha+2-n}
(E_{9,11}+E_{10,13}+E_{12,14})+
\sqrt{[\alpha+3]}q^{\alpha+3-3n}E_{15,16}).
\end{eqnarray}
\end{prop}

\end{appendix}

\end{document}